 \documentclass[12pt]{article}
  \usepackage{amsfonts}
  \usepackage[square,sort&compress,comma,numbers]{natbib} 
 \usepackage{longtable}

   \def\R{\mathbb{R}}
   \def\N{\mathbb{N}}
   \def\Z{\mathbb{Z}}
   \def\1{{\rm I\mskip -10.5mu 1}} 
   
   \def\e{{\varepsilon}}
   
   \def\lo{\mathop{\longrightarrow}}
   
   \def\vi{{\varphi}}
   \def\d{{\delta}}
   \def\l{{\lambda}}

   \def\cC{{\cal C}}

   \def\cL{{\cal L}}

   \def\supp{\mathop{\rm supp}\nolimits}
   
   \def\dist{\mathop{\rm dist}\nolimits}
   
   \def\const{\mathop{\rm const}\nolimits}
   
   \def\loc{\mathop{\rm loc}\nolimits}

   \def\no{\noindent}
   \def\proof{\mbox {{\underline {\sf Proof}} \hspace{2mm}}}
   \def\qed{{\hfill {\em q.e.d.}\\\vspace{1mm}}}
   
   \newcommand{\beq}{\begin{equation}}
   \newcommand{\eeq}{\end{equation}}

\newtheorem{df}{Definition}[section]
\newtheorem{prop}[df]{Proposition}
\newtheorem{lemma}[df]{Lemma}
\newtheorem{teo}[df]{Theorem}
\newtheorem{rem}[df]{Remark}

\newtheorem{cor}[df]{Corollary}

 \newcommand{\sezione}[1]{\section{#1}\setcounter{equation}{0}}

\begin{document}


   \title{Infinitely many positive solutions of nonlinear
     Schr\"odinger equations
}


  \maketitle


 \vspace{5mm}

\begin{center}

{ {\bf Riccardo MOLLE$^{*,a}$,\quad Donato PASSASEO$^b$}}

\vspace{5mm}

{\em
${\phantom{1}}^a$Dipartimento di Matematica,
Universit\`a di Roma ``Tor Vergata'',\linebreak
Via della Ricerca Scientifica n. 1,
00133 Roma, Italy.}

\vspace{2mm}

{\em
${\phantom{1}}^b$Dipartimento di Matematica ``E. De Giorgi'',
  Universit\`a di Lecce,\linebreak 
P.O. Box 193, 73100 Lecce, Italy.
}
\end{center}

\vspace{5mm}


{\small {\sc \noindent \ \ Abstract.} - 
\footnote{${\phantom{1}}^*$Corresponding author, {\em e-mail address:}
  molle@mat.uniroma2.it}
The paper deals with the equation $-\Delta u+a(x) u =|u|^{p-1}u 
$, $u \in H^1(\R^N)$, with  $N\ge 2$, $p>1,\ p<{N+2\over N-2}$ if
$N\ge 3$, $a\in L^{N/2}_{loc}(\R^N)$,  $\inf a>0$, $\lim_{|x| \to
  \infty} a(x)= a_\infty$. Assuming on the potential that 

 $\lim_{|x| \to \infty}[a(x)-a_\infty] e^{\eta |x|}= \infty \ \ \forall \eta>0$ and 

$ \lim_{\rho \to \infty} \sup \left\{a(\rho \theta_1) - a(\rho \theta_2)
  \ :\  \theta_1, \theta_2 \in \R^N,\  |\theta_1|= |\theta_2|=1
\right\} e^{\tilde{\eta}\rho} = 0 \hspace{0,5cm} \mbox{for some} \
\tilde{\eta}>0$, 

\noindent but not requiring any symmetry, the existence of infinitely
many positive multi-bump solutions is proved.

\noindent This result considerably improves those of  previous papers
\cite{WY,CPS1,CDS,DvS2,DPWY}.

\vspace{3mm}



\vspace{1mm}


}


\sezione{Introduction and statement of the main result}


In this paper we are concerned with the question of finding multiple
positive solutions to problem
\begin{equation}
\label{.}
\left\{
\begin{array}{l} 
-\Delta u+a(x) u  =u^p \ \ \ \mbox{in} \ \ \R^N\\
u \in H^1(\R^N),
\end{array}\right.
\end{equation}
where $N\ge 2$, $p>1,\ p< \frac{N+2}{N-2}$ if $N\ge 3$.

\no Euclidean Scalar Field equations like (\ref{.}) arise naturally in a
 large number of Physical topics like the study of solitary waves in
 nonlinear Schrodinger equations or in nonlinear Klein-Gordon
 equations. However,  besides the relevance in applied sciences, the
 interest of researchers in studying such kind of problems has been
 also due to the loss of compactness created  by the invariance of
 $\R^N$ under the action  of  translations and to the related
 challenging difficulties. 
 
 \no Actually, (\ref{.}) has a variational structure, solutions of it  can
 be searched as critical points of the functional 
\beq
\label{E}
E(u)={1\over 2}\int_{\R^N}(|Du|^2+a(x)u^2)dx-{1\over
  p+1}\int_{\R^N}|u|^{p+1}dx\qquad\forall u\in H^1(\R^N)
\eeq
 but, since $E$ does not satisfy the Palais-Smale compactness
 condition, the classical variational methods cannot be applied in a
 standard way. 
Furthermore, one can understand that the difficulty in facing problems
of this type is not only a technical fact considering that, really,
(\ref{.})  can have only the trivial solution: for instance when
the potential $a(x)$ is  increasing along a direction (see \cite{CM}).  
  
\no In this paper,  in view of their physical meaning too, we shall look
only at potentials satisfying: 
\beq
\label{1}
a\in
L^{N/2}_{\loc}(\R^N),\quad
\inf_{\R^N}a>0,\qquad\lim_{|x|\to\infty}a(x)=a_\infty.
\eeq
Starting from the Sixties of last century  many mathematicians have
devoted a lot of efforts and exploited different tools to overcome the
difficulties and to prove existence and multiplicity of solutions to
(\ref{.}).  

\no First results were obtained using the spherical symmetry of $\R^N$ and
considering radial data. 
So the existence of a ground state radial positive solution and
infinitely many radial changing sign solutions has been obtained first
by ordinary differential equations methods (see \cite{N,R}) then by
variational methods (see \cite{B,S,BeL}) taking advantage of the
compactness of the  embedding in $L^q (\R^N)$, $2<q < \frac{2N}{N-2}$,
of the subspace of $H^1(\R^N)$ consisting of radially symmetric
functions. 
It is worth also observing that under radial symmetry assumptions the
existence of infinitely many non radial changing sign solutions has
been shown (see \cite{BW}).

\no  Although analogous results could be reasonably expected when the
 symmetry in (\ref{.}) is broken by non symmetric coefficients, on the
 contrary in this case even the question of the existence appeared at
 once not easy to handle and affected by an impressive topological
 difference according the potential $a(x)$ approaches its limit at
 infinity from below or from above.   
 
\no   In the first case the existence of a positive ground state solution
  was obtained by minimizing the functional $E$ on the Nehari natural
  constraint and applying concentration-compactness arguments
  \cite{Rab,PLL}, while the multiplicity question  had an answer in
  \cite{CDS} where the existence of infinitely many changing sign
  solutions was proved assuming on the potential a decay  slower than
  any exponential decay and some stability of the directional
  derivative with respect to small perturbation of the direction. 
  
\no When $a(x)$ goes to $a_\infty$ from above the  minimization argument
does not work and, conversely, when $a(x) - a_\infty > 0$ on a
positive measure set, (\ref{.}) has not a ground state solution. 
Nevertheless, the existence of a positive bound state solution has
been shown in \cite{BL}  by subtle topological and variational
arguments, assuming a  decay of $a(x)$ faster than some exponential. 
 The multiplicity question  is even more tricky. 
 
\no   During last decade some progress has been developed  looking  mainly
  for \textit{positive  multi-bump}  solutions. 
Before discussing  nonsymmetric cases,  we mention that, again, under
symmetry assumptions the question can be controlled in a better way. 
Indeed, the existence of infinitely many positive multi-bump solutions
to (\ref{.}) has been proved assuming on $a(x)$ a suitable polynomial
decay and  \textit{radial symmetry} in \cite{WY},  \textit{planar
  symmetry} in \cite{DvS1,DPWY}. 
   
\no  The multiplicity question for (\ref{.}) involving  potentials without
 symmetry has been first considered in \cite{CPS1} where the existence
 of infinitely many positive multibump solutions (namely the existence
 for any $k\in \N$ of a $k$-bump solution) has been obtained asking to
 the potential a  ``slow''  decay with respect to some exponential plus
 a smallness of the oscillation $\sup_{x\in
   \R^N}|a-a_\infty|_{L^{N/2}(B(x,1))}$.  
However, while  a suitable decay condition on $a(x)-a_\infty$  appears
quite reasonable, the second condition seems essentially due to
technical motives. 
Hence, in subsequent papers some efforts have been made to drop this
condition, but, until now, with successful results only in the planar
case $N=2$ and assuming polynomial decay  of $a(x)$ to $a_\infty$
(\cite{DvS2,DPWY}).  
 
\no  On the other hand, it is worth remarking that a careful analysis of
 the proofs in  \cite{WY,DvS1,CPS1,DvS2,DPWY} make the reader
 understand that the symmetry in \cite{WY,DvS1}, the small oscillation
   assumption in \cite{CPS1}, the dimension restriction $N=2$ in
   \cite{DvS2,DPWY}, in spite of the different arguments and methods
   displayed in the papers, are essentially related to the same basic
   fact for the proof: working with functions having bumps located in
   regions where $a(x)- a_\infty$ is small. 
This observation is, in a way, also  validated by the results of
\cite{CMP} where the existence of infinitely many positive and
infinitely many nodal multi-bump solutions to (\ref{.}) is shown
considering potentials, having slow decay but not small oscillation
neither symmetry, which are asked to sink in some large regions of
$\R^N$  to the end of localizing the bumps suitably far and, when one
looks for changing sing solutions, to control the attractive effect of
positive and negative bumps each other. 
    
\no The result we obtain is, in our opinion, a considerable progress in
proving the existence of infinitely many positive solutions to
(\ref{.}) in non symmetric situations,  without imposing restrictions
on the dimension of the space $\R^N$ as in \cite{DvS2,DPWY} and
dropping the oscillation condition asked in \cite{CPS1}:  

\begin{teo}
\label{T} Let $a(x)$ satisfy (\ref{1})  and
\beq
\label{2}
\lim_{|x| \to \infty} [a(x) - a_\infty] e^{\eta |x|} = \infty
\hspace{1,5cm} \forall \eta > 0, \hspace{3cm} 
\eeq  
\beq
\label{1.3bis}
 \lim_{\rho \to \infty} \sup\left\{a(\rho \theta_1) - a(\rho \theta_2)
   \ :\  \theta_1, \theta_2 \in \R^N, |\theta_1|= |\theta_2|=1
 \right\} e^{\tilde{\eta}\rho} = 0 \hspace{0,5cm} \mbox{for some} \
 \tilde{\eta}>0.  \hspace{0,5cm} 
\eeq
Then problem (\ref{.}) has infinitely many solutions.
\end{teo}

\no The proof method is fully variational and it is a variant of arguments
introduced in \cite{J1,J2} and already applied in
\cite{CPS1,CMP,J3,im,CMcomp}.
Of course it is well known that solutions of (\ref{.}) correspond to
free critical points of E or, equivalently, to critical points of E on
the Nehari natural  constraint. 
However here, as in the quoted papers, critical points are searched
by min-max arguments in suitable classes of positive functions having
for all $k \in \N$ exactly $k$ ``bumps'', satisfying $k$ local Nehari
natural constraints  and $k$ local barycenter conditions. 
Therefore, being $E$ subject to constraints that are not all natural,
the min-max procedure gives rise to functions that are solutions of
equations where  Lagrange multipliers and  constraints appear and it
is an heavy task to show null the Lagrange multipliers and so proving
that  constrained critical points are actually free critical points of
$E$ and solutions of (\ref{.}). 
We point out also that, unlike the quoted papers, in the present
research the $k$-bump functions  belonging to the above described
classes must satisfy a further condition having the purpose of helping
to localize the bumps close to large radius spheres, when $k$ is
large. 

\no  Altough our  method is variational, it allows us to  describe some
 considerable asymptotic properties of the solutions we find.  
We have collected them in a proposition which, in order to be stated,
needs we introduce before the  limit problem related to (\ref{.})
\beq 
-\Delta u+a_\infty u  = |u|^{p-1}u \ \ \ \mbox{in} \ \ \R^N \ \ u \in
H^1(\R^N) 
\eeq 
and its ground state solution, denoted by $w(x)$, which is radially
symmetric, positive, unique up to translations, decreasing when the
radial coordinate increases (see f.i. \cite{BeL}). 
By the radial symmetry we shall write, with abuse of notation, $w(R)$
meaning the  value $w(x)$ takes at points $x $ such that $|x|= R$. 

\begin{prop}
\label{P1.2}
Let assumptions of Theorem \ref{T} be satisfied. 
Then $\bar{k} \in \N$ exists such that to any $ k \geq \bar{k}$ there
corresponds a positive solution $u_k$ of (\ref{.})  having the
following property: 
to $u_k$ a $k$-tuple of points $(x_1^k, \dots ,x_k^k )$  of $\R^N$
is associated in such a way that 
\beq
\label{6}
\lim_{k\to\infty}\sup\{|u_k(x+x_i^k)-w(x)|\ :\ |x|\le R, \
i=1,\ldots,k\}=0\qquad\forall R>0,
\eeq
\beq
\label{7}
\lim_{k\to\infty}\sup\{u_k(x)\ :\ |x-x_i^k|\geq R\ \mbox{ for
}i=1,\ldots,k\}=w(R)\qquad\forall R>0,
\eeq
\beq
\label{3}
\lim_{k\to\infty}\min\{|x_i^k|\ :\ i=1,\ldots,k\}=\infty,
\eeq
\beq
\label{5}
\lim_{k\to\infty}\min\{|x_i^k-x^k_j|\ :\ i\neq j,\ i,j=1,\ldots,k\}=\infty,
\eeq
\beq
\label{4}
\lim_{k\to\infty}
\frac
{\max\{|x_i^k|\ :\ i=1,\ldots,k\}}
{\min\{|x_i^k|\ :\ i=1,\ldots,k\}}
=1.
\eeq

\no Furthermore, $u_k\lo 0$ as $k\to\infty$, uniformly on the compact
subsets of $\R^N$, while $\lim_{k\to\infty}\|u_k\|_{H^1(\R^N)}$ $=\infty$ and
$\lim_{k\to\infty} E(u_k)=\infty$.
\end{prop}

\no The above proposition helps  to guess a suggestive picture of the
 solutions shape arguing in this  way:  
 the points  $(x_1^k, \dots ,x_k^k )$  are noting but the barycenters
 of the bumps which, as $k$ increases, go far from the origin and far
 away each other, while at the same time, as $k$ increases, the shape
 of $u_k$ in balls centered at $x^k_i$ for all $i=1,\dots, k$ approaches the
 shape of $w$ and outside $u_k$ decays as $w$ decays. 
So, considering the profile of $w$, one can ``see'' the $u_k$ as
functions  having an incresing number  of well glued ``bumps'' which
become more and more similar to  copies of $w$.  
Property (\ref{4}) make we understand that the bumps tend to be
distributed around spheres in $\R^N$, indeed the points ${x_1^k\over
  m_k},\ldots,{x_k^k\over m_k}$, where $m_k=\min\{|x^k_i| $:
$i=1,\ldots,k\}$,  as $k\to \infty$ become closer and closer to the
sphere of radius 1 centered at the origin  $ \partial
B(0,1)\subset\R^N$. 
Furthermore, as we shall see in Corollary \ref{Cnew}, more  can be
asserted, namely that the  distribution of 
these points tends to be ``uniform'' around $ \partial B(0,1)$,
because, for all $x\in \partial B(0,1) $ and for all $r>0$ the number
of the points ${x_1^k\over m_k},\ldots,{x_k^k\over m_k}$ lying in
$B(x,r)$ tends to infinity, as $k$ goes to infinity, with rate $k\,
r^{N-1}$.  

\no The paper is organized as follows. 
In Section 2 the classes of
$k$-bumps functions in which the solutions are seeked are introduced
and their properties are recalled. 
In Section 3 the min-max arguments to find the good candidates to  be
critical points are displayed, and in Section 4 the asymptotic
behaviour of these functions is described as the number of the bumps
increses. 
Finally, in Section 5  the before found $k$-bump functions are shown to
be free critical points of $E$. 

\vspace{2mm}

{\sc Acknowledgements.} The Authors would like to thank very much
Prof. Giovanna Cerami for many useful discussion on this subject.


\sezione{Variational framework and known facts}

    
 Throughout the paper we make use of the following notation:
 
 \begin{itemize}
\item $H^{1}({\mathbb R}^{N})$ is the usual Sobolev space endowed with
  the standard scalar product and norm 
\begin{displaymath}
(u, v)=\int_{\mathbb R^N}[\nabla u \nabla v+a_\infty uv]dx;\qquad
\|u\|^{2}=\int_{{\mathbb R}^N}\left[|\nabla u|^{2}+a_\infty u^{2}\right]dx; 
\end{displaymath}
\item $L^q(\Omega)$, $1\leq q \leq +\infty$, $\Omega \subseteq \mathbb
  R^N$, denotes a Lebesgue space, the norm in $L^q(\Omega)$ is denoted
  by $|u|_{q, \Omega}$ when $\Omega$ is a proper subset of $\mathbb
  R^N$, by $|\cdot|_q$ when $\Omega=\mathbb R^N$; 
\item for any $\rho>0$ and for any $z\in \mathbb R^N$, $B(z,\rho)$
  denotes the ball of radius $\rho$ centered at $z,$ and
  $S(z,\rho)=\partial B(z,\rho)$; 
\item for any measurable set $ \mathcal{O} \subset \R^N, \
  |\mathcal{O}|$ denotes its Lebesgue measure; 
\item $c,c', C, C', C_i,\ldots$ denote various positive constants.
    
\end{itemize}

\no Moreover, in Appendix we report a table of the main notations we
use in this paper. 
   
\no In what follows we denote by 
\beq
E_\infty(u)={1\over 2}\int_{\R^N}(|D u|^2+a_\infty u^2)dx-{1\over
  p+1}\int_{\R^N}|u|^{p+1}dx,\qquad u\in H^1(\R^N),
\eeq
 the functional related to the limit problem. 
Next lemma (see \cite{CMilan,AM} and the references
therein) summarizes the main properties of the ground state solution
$w$ of it.
\begin{lemma}
\label{Lw}
The function $w$ is unique up to translations, has radial symmetry,
decreases when the radial 
coordinate increases and satisfies
\begin{equation}
\lim_{|x|\to\infty}|D^i w
(x)||x|^{\frac{N-1}{2}}e^{\sqrt{a_\infty}|x|}=d_i>0\quad \mbox{ with }
 d_i\in\R\ \mbox{ for } i=0,1.
\end{equation}
If we set
\beq
\label{2.2}
w_y(x)=w(x-y)\quad\forall x ,y\in\R^N\ \mbox{ and }Z=\{w_y\ :\
y\in\R^N\},
\eeq
then $Z$ is non degenerate, namely the following properties are true: 
\begin{itemize}
\item[{\em a)}] $(E_{\infty})''(w_y)$ is an index zero Fredholm map
  for all $y \in \R^N$; 
\item[{\em b)}] Ker $(E_{\infty})''(w_y)$ = span $\left\{
    \frac{\partial w_y}{\partial x_j}:  j=1,\ldots,N \right\} =
  T_{w_y}(Z)$, where $T_{w_y}(Z)$ is the tangent space to $Z$ at
  $w_y$.
\end{itemize}
\end{lemma}

\no Let us set $a_0=\inf_{\R^N}a$ and fix $\d>0$ such that
\beq
\label{delta}
\delta<\min\left\{
1,
\left( a_0\over p\right)^{1\over p-1},
\left( a_0\over 2\right)^{1\over p-1},
{a_0^{1\over p-1}\over 2},
{w(0)\over 3}\right\},
\eeq
then we denote by $R_\d>0$ a fixed number so that $w(x)<\d$ $\forall x\in
\R^N\setminus B(0, R_\d/2)$.

\no For every function $u\in H^1(\R^N)$, $u\ge 0$, we denote by
\beq
u_\delta=u\wedge\d,\qquad u^\d=u-u_\d
\eeq
\no 
and call $u_\d$ and $u^\d$ the submerged and the emerging part of $u$.
We say that a function $u\in
H^1(\R^N)$ is emerging  around $x_1,\ldots,x_k\in \R^N$ if 
$u^\d=\sum_{i=1}^ku^\d_i$ where, for all $i\in\{1,\ldots,k\}$,
$u^\d_i(x)=0$ $\forall x\not\in B(x_i,R_\d)$ and
$u^\d_i\not\equiv 0$.
On the submerged parts, the functional $E$ has the following features.

\begin{rem}
\label{C}
{\em The functional $E$ is coercive and convex, hence weakly lower semicontinuous,
on the convex set 
\beq
\cC=\{u\in H^1(\R^N)\ :\ |u|\le\delta\ \mbox{ a.e. in }\R^N\}.
\eeq
}\end{rem}

\no Indeed, taking into account the choice of $\d$ in (\ref{delta}), we
have
\begin{eqnarray}
\nonumber E(u) & \ge & {1\over 2} \int_{\R^N} (|D u|^2+a_0 u^2)dx-{1\over
  p+1}\int_{\R^N}|u|^{p+1}dx \\
 &\ge & 
 {1\over 2} \int_{\R^N} |D u|^2 dx +\left({a_0\over 2}-{\d^{p-1}\over
     p+1}\right) \int_{\R^N}  u^2 dx\\
\nonumber &\ge &  c\, \|u\|^2\qquad\forall u\in\cC
\end{eqnarray}
for $c>0$, and
\beq
{d\over dt}E(u_1+t(u_2-u_1))\ge \int_{\R^N} \hspace{-3mm} |D (u_2-u_1)|^2dx
+\int_{\R^N}\hspace{-3mm} (a_0-p\d^{p-1})(u_2-u_1)^2dx\ge 0\quad\forall
u_1,u_2\in\cC.
\eeq

\vspace{2mm}

\no For all $k\ge 1$, we  set 
\beq
\label{2.3}
\left\{
\begin{array}{lc}
D_1=\R^N & \mbox{if }k=1\\
D_k=\{(x_1,\ldots,x_k)\in(\R^N)^k\ :\ |x_i-x_j|\ge 3R_\d\ \mbox{ for
}i\neq j,\ i,j=1,\ldots k\}& \mbox{if }k\ge 2
\end{array}
\right.\eeq
and for all $(x_1,\ldots,x_k)\in D_k$ we consider the set consisting
of functions emerging around $x_1,\ldots,x_k$ and satisfying  local
Nehari and local barycenter constraints:
\begin{eqnarray}
\nonumber S_{x_1,\ldots,x_k }\hspace{-2mm} &=\{&\hspace{-2mm}
u\in H^1(\R^N)\  : u \geq  0, \  u \mbox{
   emerging around }x_1,\ldots,x_k,\\
& &
 \label{2.6} E'(u)[u_i^\delta]=0,\  \beta_i(u)=x_i,\ \forall i\in\{1,\ldots,k\}\},
\end{eqnarray}
where $\beta_i$ is the local barycenter defined by
\beq
\label{2.4}
\beta_i(u)={1\over |u^\d_i|_2^2} \int_{\R^N}x\,
(u^\d_i(x))^2\,dx\qquad \mbox{ for }i=1,\ldots,k.
\eeq.

\no Analogously, we denote by $S^\infty_y$, $y\in\R^N$, the set
obtained replacing $E$ by $E_\infty$ in (\ref{2.6}).

\no Notice that, if $u\in S_{x_1,\ldots,x_k }$, it satisfies the equality
\beq
\label{2.5bis}
\hspace{-2cm}
\beta'_i(u)[\psi]=
2\left(\int_{\R^N}(u^\d_i(x))^2\, dx\right)^{-1}
\int_{\R^N} u^\d_i(x)\psi(x)\, (x-x_i)\, dx
\eeq
$$
\qquad \hspace{3cm} \forall \psi\in H^1_0(B(x_i,R_\d)),\ \forall
i\in\{1,\ldots,k\}, 
$$
as one can verify by direct computation.

\vspace{2mm}

\no In the following statements we collect some features, whose proof can
be found in \cite{CMP} (see also \cite{CPS1}), that draw the
variational setting we are working in.

\begin{prop} {\em (see \cite[Lemma 2.10]{CMP})}
\label{Pproiezioni} 
Let $u\in H^1(\R^N)$ be such that $u^\d\neq 0$ and $u^\d$ has compact
support.
Then there exists a unique $\bar t \in (0,\infty)$ such that
$E'(u_\d+\bar t u^\d)[u^\d]=0$; moreover $\bar t$ is the maximum point
of the function $t\mapsto E(u_\d+tu^\d)$.

\no The same statements hold when we consider the functional $E_\infty$.
\end{prop}  

\no As a consequence of Proposition \ref{Pproiezioni}, it is readily
seen that $S_{x_1,\ldots,x_k}\neq\emptyset$ for every
$(x_1,\ldots,x_k)\in D_k$. 
Indeed, let $\phi\in\cC^\infty_0(B(0,R_\d))$ be a positive radially
symmetric function such that $\phi^\d\neq 0$, then 
$\sum_{i=1}^k[(\phi(x-x_i))_\d+\bar t_i(\phi(x-x_i))^\d]\in
S_{x_1,\ldots,x_k}$, where $\bar t_i$ is the value corresponding to 
$\phi(x-x_i)$ provided by Proposition  \ref{Pproiezioni}.

\vspace{1mm}

\no  For every function $v\in H^1(\R^N)$ with compact support, let
us define
\begin{eqnarray}
\nonumber F(v)&=&\frac{1}{2}\int_{\R^N}(|D v|^2 +a(x)v^2)dx + 
\int_{ \R^N} a(x) \delta \, v \, dx  \hspace{6ex} 
 \\
& &- \frac{1}{p+1}\int_{\supp\, v} (\delta + v)^{p+1} dx +
\frac{\delta ^{p+1}}{p+1}\  |\supp v|.  
\end{eqnarray}
By the choice of $\d$, for every $(x_1,\ldots,x_k)\in D_k$ we can write
\beq
E(u)=E(u_\d)+F(u^\d)
=E(u_\d)+\sum_{i=1}^k F(u^\d_i)\qquad\forall u\in S_{x_1,\ldots,x_k}
\eeq
and both $E$ and $F$ have positive sign on our set of functions:

\begin{prop} {\em (see \cite[Proposition 2.15]{CMP})} 
\label{P2.15CMP}
Assume $(x_1,\ldots,x_k)\in D_k$ and  $u\in S_{x_1,\ldots,x_k}$, then
\beq
E(u)>0,\qquad F(u^\d_i)>0\quad \forall i\in\{1,\ldots,k\}.
\eeq
\end{prop}

\vspace{2mm}

\begin{df}
{\em
For every $(x_1,\ldots,x_k)\in D_k$ we set 
\beq
\label{df}
f_k(x_1,\ldots,x_k):=\inf\{E(u)\ :\ u\in S_{x_1,\ldots,x_k}\}.
\eeq
}
\end{df}

\begin{prop}
\label{Pnew}
For every  $(x_1,\ldots,x_k)\in D_k$ the infimum in (\ref{df}) is achieved and 
\beq
f_k(x_1,\ldots,x_k)>0.
\eeq
Moreover, for every $\bar u\in S_{x_1,\ldots,x_k}$ such that $E(\bar
u)=f_k(x_1,\ldots,x_k)$ the following properties hold:
\begin{itemize}
\item[i)]
$\bar u(x)>0$ $\forall x\in\R^N$ and $\bar u(x)<\d$ for all $x$ such that $\dist(x,\supp \bar u^\d)>0$;
\item[ii)] $\bar u$ satisfies the equation 
\beq
\label{e}
-\Delta\bar u(x)+a(x)\bar u(x)=\bar u^p(x)
\qquad \forall x\in \R^N\ \mbox{ s.t. }\dist(x,\supp \bar u^\d)>0;
\eeq
\item[iii)]
there exist two positive constants $b$ and $c$ such that 
\beq
\label{dec}
\bar u(x)\le c\, e^{-b\, d(x)}\qquad\forall x\in
\R^N\setminus\supp\bar u^\d
\eeq
where $d(x)=\dist (x,\supp\bar u^\d)$;
\item[iv)]
there exist Lagrange multipliers
$\lambda_1,\ldots,\lambda_k$ in $\R^N$ such that  
\beq
\label{lambda}
E'(\bar u)[\psi]=\int_{B(x_i,R_\d)}\hspace{-4mm}\bar u^\d_i
(x)\psi(x)[\lambda_i\cdot (x-x_i)]dx\quad \forall \psi\in
H^1_0(B(x_i,R_\d)),\ \forall i\in\{1,\ldots,k\}. 
\eeq
\end{itemize}
\end{prop}

\no The existence of a minimizer $\bar u\in S_{x_1,\ldots,x_k}$ is
proved in \cite[Proposition 3.1]{CMP}, {\em (i) -- (iii)} are
contained in \cite[Lemma 3.4]{CMP} (for the property $\bar u>0$ see also
\cite[Lemma 3.4]{CPS1}) and {\em (iv)} is in \cite[Proposition
3.5]{CMP}.

\begin{rem}
{\em
 For every $R\in \left[R_\d,{1\over 2}\min\{|x_i-x_j|\ :\ i\neq j,\
i,j=1,\ldots,k\}\right]$, the relation
\beq
\label{B}
\sup\{\bar u(x)\ :\ x\in\R^N\setminus\cup_{i=1}^k B(x_i,R)\}
=\sup\{\bar u(x)\ :\ x\in\partial(\cup_{i=1}^k B(x_i,R))\}
\eeq
holds true.
Indeed, since $u$ is a positive solution of (\ref{e})
in $\R^N\setminus\supp (u^\d)$, considering the choice of  $\d$ we get
$\Delta\bar u>0$  in $\R^N\setminus\cup_{i=1}^k B(x_i,R)$.
So maximum principle gives (\ref{B}). 
}\end{rem}

\no Concerning the limit problem, we have: 

\begin{lemma} {\em (see \cite[Lemma 4.1]{CMP})}
\label{G}
For every $y\in\R^N$, $w_y\in S^\infty_y$ (see (\ref{2.2})) and
\beq
\min_{S^\infty_y}E_\infty =E_\infty(w_y).
\eeq
\end{lemma}

\no Now, we describe the asymptotic behaviour of suitable
sequences of minimizing functions and of the corresponding Lagrange
multipliers. 

\begin{prop}
\label{Pw}
Let $(k_n)_n$ be a sequence in $\N$ and
$((x_{1,n},\ldots,x_{k_n,n}))_{n}$ a sequence such that   
\beq
\label{bo}
\lim_{n\to \infty}\min\{|x_{i,n}|\ :\ i=1,\ldots,k_n\}=\infty,
\eeq
\beq
\lim_{n\to \infty}\min\{|x_{i,n}-x_{j,n}|\ :\ i\neq j,\
i,j=1,\ldots,k_n\}=\infty.
\eeq
For all $n\in\N$, let $u_n$ be a minimizer of $E$ in
$S_{x_{1,n},\ldots,x_{k_n,n} }$. 
Let $\lambda_{i,n}$ be the Lagrange multipliers provided by (iv) of
Proposition \ref{Pnew}, then 
\beq
\label{w}
\lim_{n\to\infty}
\sup\{|u_n(x+x_{i,n})-w(x)|\ :\ |x|\le R,\
i=1,\ldots,k_n\}=0\qquad\forall R>0,
\eeq
\beq
\label{1527}
u_n(x+x_{i,n})\to w(x)\qquad\mbox{ in } H^1_{\loc}(\R^N)
\eeq
 and
\beq
\label{lambdan}
\lim_{n\to \infty}\max\{|\lambda_{i,n}|\ :\ i=1,\ldots,k_n\}=0.
\eeq
\end{prop}

\no For the proof  we refer the reader to \cite[Proposition 5.5]{CMP}.
In fact, (\ref{1527}) is {\em(b)} in the proof of Proposition 5.5 in
\cite{CMP} while (\ref{lambdan}) here corresponds to (5.32) in
that proof.

\vspace{1mm}

\no Finally, let us prove the following continuity property.

\begin{prop}
\label{P1627}
Let $a(x)$ verify assumptions (\ref{1}), then $f_k:D_k\to \R$  is a
continuous function. 
\end{prop}

\no \proof 
The upper semicontinuity is proved in \cite[Lemma 4.2]{CMP}.

\no In order to prove the lower semicontinuity, let us consider a
sequence $((x_1^n,\ldots,x_k^n))_n$ in $D_k$  such that 
$(x_1^n,\ldots,x_k^n)\to(x_1,\ldots,x_k)\in D_k$, as $n\to\infty$, and
let $u_n\in S_{x_1^n,\ldots,x_k^n} $ be such that
$E(u_n)=f_k(x_1^n,\ldots,x_k^n)$.

\no Since $f_k$ is upper semicontinuous, $E(u_n)$ is bounded and so we
can infer that $F(u_n^\d)$ is bounded and $((u_n)_\d)_n$ is bounded in
$H^1_0(\R^N)$, taking into account $E(u_n)=E((u_n)_\d)+F(u_n^\d)$,
Proposition \ref{P2.15CMP} and the coercivity of $E$ on the submerged
parts (see Remark \ref{C}).  

\no
Let us show that also $(u_n^\d)_n$ is bounded in $H^1(\R^N)$, that is
$((u^\d_i)_n)_n$ is bounded for every $i\in\{1,\ldots,k\}$. 
From $E'(u_n)[(u^\d_i)_n]= 0$, $\forall n\in\N$,  we get
\beq
\label{1771}
\int_{\R^N}(|D (u^\d_i)_n|^2+a(x)(u^\d_i)_n^2)dx+
\delta\int_{\R^N}a(x)(u^\d_i)_ndx=\int_{\R^N}(\d+(u^\d_i)_n)^p(u^\d_i)_ndx.
\eeq
Hence we can write
\beq
\label{1741}
\begin{array}{rcl}
\vspace{2mm}
F((u^\d_i)_n)&=&{1\over 2} \int_{\R^N}(\d+(u^\d_i)_n)^p(u^\d_i)_ndx
+{1\over 2} \int_{\R^N}a(x)(u^\d_i)_ndx\\
& & -{1\over p+1} \int_{\supp
  (u^\d_i)_n}(\d+(u^\d_i)_n)^{p+1}dx+{\d^{p+1}\over p+1} \, |\supp 
(u^\d_i)_n|  
\end{array}
\eeq
and taking into account that  $|\supp
(u^\d_i)_n|\le c_1$, $\forall n\in\N$,  we see that
\beq
\label{1740}
\begin{array}{rcl}
\vspace{2mm}
F((u^\d_i)_n)& \ge & {1\over 2} \int_{\R^N}(\d+(u^\d_i)_n)^p(u^\d_i)_ndx
-{1\over p+1} \int_{\supp (u^\d_i)_n}(\d+(u^\d_i)_n)^{p+1}dx\\
\vspace{2mm}
& \ge &\int_{\{(u^\d_i)_n>{4\d\over
    p-1}\}}(\d+(u^\d_i)_n)^p(u^\d_i)_n\left[{1\over 
    2}-{1\over p+1}\left({\d\over (u^\d_i)_n}+1\right)\right]dx\\
\vspace{2mm}
& & -{1\over p+1}\int_{\{0<(u^\d_i)_n\le {4\d\over
    p-1}\}}(\d+(u^\d_i)_n)^{p+1}dx\\
& \ge &\int_{\{(u^\d_i)_n>{4\d\over p-1}\}}((u^\d_i)_n)^{p+1}\left[{1\over
    4}\, {p-1\over p+1}\right]dx-c_2.
\end{array}
\eeq
From $F((u^\d_i)_n)\le c_3$ and  (\ref{1740}) it follows 
\beq
\int_{\{(u^\d_i)_n>{4\d\over p-1}\}}((u^\d_i)_n)^{p+1}dx\le\const,
\eeq
so that $((u^\d_i)_n)_n$ is bounded in $L^{p+1}$, in $L^2$, in $L^1$
and so it turns out to be bounded also in $H^1$ by (\ref{1771}).

\no Summarizing, $(u_n)_n$ is bounded in $H^1$ so, up to a subsequence,
it converges to a function $\bar u$ weakly in $H^1(\R^N)$ and we have
also that $(u^\d_i)_n\to  \bar u_i^\d$ strongly in $L^{p+1}$ and in $L^{2}$ .

\no Observe that  $\bar u^\d_i\neq 0$, indeed from (\ref{1771}) and
the choice of $\delta$ in (\ref{delta}) we obtain 
\beq
\begin{array}{rcl}
\vspace{2mm}
0 & \ge & c_1|(u^\d_i)_n|_{p+1}^2+(a_0\d-2^{p-1}\d^p)\int_{\supp
  (u^\d_i)_n}(u^\d_i)_n dx-2^{p-1}|(u^\d_i)_n|^{p+1}_{p+1} \\
& \ge & c_1|(u^\d_i)_n|_{p+1}^2-c_2|(u^\d_i)_n|_{p+1}^{p+1}
\end{array}
\eeq
that implies $|(u^\d_i)_n|_{p+1}\ge\const$ $\forall n\in\N$.
Then, $\bar u$ is a function emerging around
$(x_1,\ldots,x_k)$ and it verifies $\beta_i(\bar u)=x_i$ $\forall
i\in\{1,\ldots,k\}$, by the $L^2$-convergence of the emerging parts. 

\no Now, according to Proposition \ref{Pproiezioni}, let
$t_i\in(0,\infty)$, $i\in\{1,\ldots,k\}$, be such that $\hat
u=\bar u_\d+\sum_{i=1}^kt_i\bar u^\d_i \in S_{x_1,\ldots,x_k}$.
We have
\beq
\begin{array}{rcl}
\vspace{2mm}
f(x_1,\ldots,x_k) & \le & E(\hat u)=E(\bar
u_\d+\sum_{i=1}^kt_i\bar u^\d_i)\\
\vspace{2mm}
& \le & \liminf\limits_{n\to\infty}E( (u_\d)_n+ \sum_{i=1}^kt_i(u^\d_i)_n)
\\
\vspace{2mm}
& \le &  \liminf\limits_{n\to\infty}E( u_n) =
\liminf\limits_{n\to\infty}f(x_1^n,\ldots,x_k^n), 
\end{array}
\eeq
that is the desired conclusion.

\qed


\sezione{The min-max argument}


\no In this section we use a min-max argument as in \cite{Lincei} to
 obtain suitable $k$-bumps functions that in next section will be
 proved to be solutions.

 For all $\sigma>0$  and for all $k\ge
2$, let us set
\beq
\label{2.17}
\begin{array}{rcl}
\vspace{2mm}
D^{k,\sigma}(\rho,\theta_1,\ldots,\theta_k)&=&\left\{(x_1,\ldots,x_k)\in D_k\
:\ \left[{1\over k}\sum_{i=1}^k|x_i|^2\right]^{1/2}=\rho,\right.
\\ \vspace{2mm}
& &
\left.
(1+2\sigma)^{-1}\rho\le |x_i|\le (1+2\sigma)\rho,\ 
 {x_i\over |x_i|}=\theta_i\ \mbox{ for }i=1,\ldots,k\right\}\
\\ 
& &
\hspace{3cm}\forall\rho\ge 0,\ \forall (\theta_1,\ldots,\theta_k)\in[S(0,1)]^k
\end{array}
\eeq
and 
\beq
\label{2.18}
D^{k,\sigma}=  \{(\rho,\theta_1,\ldots,\theta_k)\ :\ \rho\ge 0,\
\theta_i\in S(0,1)\ \mbox{ for }i=1,\ldots,k,\
D^{k,\sigma}(\rho,\theta_1,\ldots,\theta_k)\neq\emptyset\}.
\eeq
Then, let us consider the continuous function
$g^{k,\sigma}:D^{k,\sigma}\to \R$ defined by
\beq
\label{g}
g^{k,\sigma}(\rho,\theta_1,\ldots,\theta_k)=\min\{f_k(x_1,\ldots,x_k)\
:\ (x_1,\ldots,x_k)\in D^{k,\sigma}(\rho,\theta_1,\ldots,\theta_k)\},
\eeq
where the minimum is achieved because $f_k$ is a continuous
function and $D^{k,\sigma}(\rho,$ $\theta_1,\ldots,$ $\theta_k)$ is a
compact subset of $(\R^N)^k$.
The number $\sigma>0$ will be fixed later.

\begin{prop}
\label{P2.2}
Assume that the potential $a(x)$ satisfies conditions (\ref{1}) and
(\ref{2}) and let $k\ge 2$.
Then,
\beq
\label{1120}
\sup_{D^{k,\sigma}}g^{k,\sigma}>k\, E_\infty(w).
\eeq
Moreover, there exist $r_k>0$ and $(x_1^k,\ldots,
x_k^k)$ in $D_k$ such that 
\beq
\label{2.20}
\left[{1\over k}\sum_{i=1}^k|x_i^k|^2\right]^{1/2}=r_k,\quad 
(1+2\sigma)^{-1}r_k\le |x_i^k|\le (1+2\sigma) r_k,
\mbox{ for }i=1,\ldots,k
\eeq
and
\beq
\label{2.21}
f_k(x_1^k,\ldots, x_k^k)=g^{k,\sigma}\left(r_k,
{x_1^k\over |x_1^k |},\ldots, {x_k^k\over|
  x_k^k|}\right)=\max_{D^{k,\sigma}}g^{k,\sigma}. 
\eeq
\end{prop}

\no\proof
In order to prove that (\ref{1120}) holds, let us choose
$\tilde\theta_1,\ldots,\tilde\theta_k$ 
  in $S(0,1)$ such that $\tilde\theta_i\neq\tilde\theta_j$
  for $i\neq j$.
Then, there exists $\tilde\rho>0$ such that $(\rho,
\tilde\theta_1,\ldots,\tilde\theta_k)\in D^{k,\sigma}$ $\forall
\rho\ge\tilde\rho$.
For all $\rho\ge\tilde\rho$, choose $(\tilde x_{1,\rho},\ldots,\tilde
x_{k,\rho})$ in
$D^{k,\sigma}(\rho,\tilde\theta_1,\ldots,\tilde\theta_k)$ and $\tilde
u_\rho$ in $S_{\tilde x_{1,\rho},\ldots,\tilde x_{k,\rho}}$ such that 
\beq
\label{e3.7}
E(\tilde
u_\rho)=f_k(\tilde x_{1,\rho},\ldots,\tilde
x_{k,\rho})=g^{k,\sigma}(\rho,\tilde\theta_1,\ldots,\tilde\theta_k). 
\eeq
Notice that $\lim_{\rho\to\infty}|\tilde x_{i,\rho}|=\infty$ for
$i=1,\ldots,k$ and $\lim_{\rho\to\infty}|\tilde x_{i,\rho}-\tilde
x_{j,\rho}|=\infty$ for $i\neq j$ because $\tilde
\theta_i\neq\tilde\theta_j$.
Then, since $a(x)\to a_\infty$, by (\ref{1527}),
(\ref{dec}) and (\ref{e}), we obtain
\beq
\label{kEinfty}
\begin{array}{rcl}
\vspace{3mm}
\lim\limits_{\rho\to\infty}  E(\tilde u_\rho)&=&\lim\limits_{\rho\to\infty}
\sum_{i=1}^k E(w(\cdot-\tilde x_{i,\rho})) \\
&=&\lim\limits_{\rho\to\infty} \sum_{i=1}^k E_\infty (w(\cdot-\tilde
x_{i,\rho}))=k E_\infty (w). 
\end{array}
\eeq

\no Our next goal is to show that $E(\tilde u_\rho)$ approaches $k E_\infty
(w)$ from above as $\rho\to\infty$.
Notice that 
\beq
\liminf_{\rho\to\infty}{|\tilde x_{i,\rho}-\tilde
x_{j,\rho}|\over \rho} >0\qquad\mbox{ for }i\neq j
\eeq
because $\tilde \theta_i\neq\tilde\theta_j$ and that, if we set 
\beq
r_\rho={1\over 4}\min\big\{|\tilde x_{i,\rho}-\tilde x_{j,\rho}|\ :\
  i,j\in\{1,\ldots,k\},\ i\neq j\big\},
\eeq
we have
\beq
\label{2.27}
\liminf_{\rho\to\infty}{r_\rho\over\rho}>0.
\eeq
Then, let us define 
\beq
u_{\rho,i}(x)=\zeta\left({|x-\tilde x_{i,\rho}|\over r_\rho}\right)\,
\tilde u_\rho(x),\qquad x\in\R^N,
\eeq
where $\zeta\in\cC^\infty([0,\infty),[0,1])$ is a cut-off function
such that $\zeta(t)=1$ if $t\in[0,1]$, $\zeta(t)=0$ if
$t\in[2,\infty)$.
By (\ref{dec}), (\ref{e}) and (\ref{2.27}) there exist constants
$c,c_k>0$ such that 
\beq
\label{1419}
E(u_\rho)=\sum_{i=1}^k
E(u_{\rho,i})+O(e^{-c\, r_\rho})=\sum_{i=1}^kE(u_{\rho,i})+O(e^{-c_k\rho}),
\eeq
where the constant $c_k$ depends only on  $\tilde
\theta_1,\ldots,\tilde\theta_k$ and $\sigma$.
In order to evaluate $E(u_{\rho,i})$, let us consider
$t^\infty_i\in (0,\infty)$ such that
$v_{\rho,i}:=(u_{\rho,i})_\d+t^\infty_i u^\d_{\rho,i}\in
S^\infty_{\tilde x_{i,\rho}}$ (see Proposition \ref{Pproiezioni}).
Then, taking also into account Proposition \ref{Pproiezioni} and Lemma
\ref{G}, for every $i\in\{1,\ldots,k\}$ and large $\rho$ we have
\beq
\label{1420}
\begin{array}{rcl}
E(u_{\rho,i})& \ge&  E(v_{\rho,i})  = 
E_\infty (v_{\rho,i})+\displaystyle{\int_{B(\tilde
  x_{i,\rho},2r_\rho)}(a(x)-a_\infty) v_{i,\rho}^2dx}
\\[3ex]
&\ge & 
E_\infty (w)+\displaystyle{\int_{B(\tilde
  x_{\rho,i},R_\d)}(a(x)-a_\infty) (v_{\rho,i}^\delta)^2dx}.
\end{array}
\eeq
By (\ref{1527}), $|v^\d_{\rho,i}-w^\d_{\tilde x_{i,\rho}}|_{2, B(\tilde
  x_{\rho,i},R_\d)}\to
0$ as $\rho\to\infty$, then (\ref{2}) implies
\beq
\label{1421}
\lim_{\rho\to\infty}\left[\int_{B(\tilde
    x_{i,\rho},R_\d)}(a(x)-a_\infty)(v^\d_{\rho,i})^2dx\right]e^{c_k\rho}=\infty. 
\eeq
Setting $\alpha(\rho)=\sum_{i=1}^k\int_{B(\tilde
    x_{i,\rho},R_\d)}(a(x)-a_\infty)(v^\d_{\rho,i})^2dx$, from
  (\ref{1419}) and (\ref{1420}) it follows
\beq
E(u_\rho)\ge k\, E_\infty(w)+\alpha(\rho)+O(e^{-c_k\rho}),
\eeq
so we have the desired asymptotic behaviour and (\ref{1120}) follows
from (\ref{1421}).

\no Now, in order to prove that $\max_{D^{k,\sigma}} g^{k,\sigma}$ is
  achieved, consider a sequence
  $((\rho_n,\theta_{1,n},\ldots,$ $\theta_{k,n}))_n$ in
  $D^{k,\sigma}$ such that
\beq
\lim_{n\to\infty} g^{k,\sigma}(
\rho_n,\theta_{1,n},\ldots,\theta_{k,n})=\sup_{D^{k,\sigma}}
g^{k,\sigma}.
\eeq
Let us prove that the sequence $(\rho_n)_n$ is bounded.

\no Arguing by contradiction, assume that, up to a subsequence,
$\lim_{n\to\infty}\rho_n=\infty$ and, for all $n\in\N$, choose
\beq
(x_{1,n},\ldots, x_{k,n})\in D^{k,\sigma}
(\rho_n,\theta_{1,n},\ldots,\theta_{k,n})
\eeq
such that $f_k(x_{1,n},\ldots, x_{k,n} )=g^{k,\sigma}(
\rho_n,\theta_{1,n},\ldots,\theta_{k,n})$.

\no For every $i\in\{1,\ldots,k\}$ and $n\in\N$, let $t_{i,n}\in
(0,\infty)$ be such that $\tilde w_{x_{i,n}}:=(w_{x_{i,n}})_\d+t_{i,n}
w_{x_{i,n}}^\d\in S_{x_{i,n}}$.
Notice that since $\rho_n\to\infty$, as $n\to\infty$, and $a(x)\lo
a_\infty$, as $|x|\to\infty$, then $t_{i,n}\to 1$, by its definition, so that
$\|\tilde w_{x_{i,n}}-w_{x_{i,n}}\|\to 0$.
Hence
\beq
\label{1114}
\lim_{n\to\infty}E(\tilde w_{x_{i,n}} )=E(w).
\eeq
We observe that
\beq
\tilde w_{x_{1,n}}\vee\ldots\vee
\tilde w_{x_{k,n}}\in S_{x_{1,n},\ldots,x_{k,n} },
\eeq
and so, by the coercivity of $E$ on the submerged parts, we obtain
\beq
\label{117}
\begin{array}{rcl}
f_k(x_{1,n},\ldots,x_{k,n} ) &\le &
E(\tilde w_{x_{1,n}}\vee\ldots\vee
\tilde w_{x_{k,n}} )\\[1ex]
& = & E(\tilde w_{x_{1,n}})+E (\tilde w_{x_{2,n}}\vee\ldots\vee
\tilde w_{x_{k,n}} )\\[1ex]
& & -E(\tilde w_{x_{1,n}}\wedge (\tilde w_{x_{2,n}}\vee\ldots\vee
\tilde w_{x_{k,n}}) )\\[1ex]
&\le &E(\tilde w_{x_{1,n}})+E (\tilde w_{x_{2,n}}\vee\ldots\vee
\tilde w_{x_{k,n}} )\\[1ex]
\vdots \phantom{****}&\vdots&\phantom{****} \vdots\\[1ex]
&\le &\sum_{i=1}^k E(\tilde w_{x_{i,n}} ).
\end{array}
\eeq
By (\ref{1114}) and (\ref{117}) we infer that
\beq
\limsup_{n\to\infty} f_k(x_{1,n},\ldots, x_{k,n} )\le k\, E_\infty(w),
\eeq
which is in contradiction with (\ref{1120}).
Therefore, the sequence $(\rho_n)_n$ must be bounded and (up to a
subsequence)
\beq
\lim_{n\to\infty}\rho_n=r_k,\quad 
\lim_{n\to\infty}(x_{1,n},\ldots,x_{k,n})=(x^k_1,\ldots,x^k_k)
\eeq
for suitable $r_k>0$ and $(x^k_1,\ldots,x^k_k)$ in $D_k$.
Thus, all the assertions of Proposition \ref{P2.2} hold for $r_k$
and $(x^k_1,\ldots,x^k_k)$, by the continuity of $f_k$ (see Proposition \ref{P1627}).

\qed

\vspace{2mm}

\no Our aim will be to prove that every function $u_k\in
S_{x^k_1,\ldots,x^k_k}$,  such that
$E(u_k)=f_k(x^k_1,\ldots,$ $x^k_k)=\max_{D^{k,\sigma}}g^{k,\sigma}$,
is a solution of problem (\ref{.}) for $k$ large enough and $\sigma>0$
suitably chosen.

\begin{rem}
\label{R2.3a}
{\em
Notice that in the proof of Proposition \ref{P2.2} the existence of
the positive constant $c_k$ in (\ref{1419}) is strictly related to the fact
that the minimum 
\beq
\tilde\mu(\tilde\theta_1,\ldots,\tilde\theta_k)=\min\{|\tilde\theta_i-\tilde\theta_j|\
:\ i,j\in\{1,\ldots,k\},\ i\neq j\}
\eeq
is positive. Moreover, as one can verify by direct computation,
$ c_k\to 0$ as
$\tilde\mu(\tilde\theta_1,\ldots,\tilde\theta_k)$ $\to 0$.

\no On the other hand, it is clear that the maximum 
\beq
\tilde\mu_k=\max\{\tilde\mu(\tilde\theta_1,\ldots,\tilde\theta_k)\ :\
\tilde\theta_i\in S(0,1)\ \mbox{ for }i=1,\ldots,k\}
\eeq
tends to 0 as $k\to\infty$.
Therefore, $ c_k$ must tend to 0 as $k\to\infty$.
This fact explains why in this paper  we need
condition (\ref{2}), while the decay condition 
\beq
\exists \ \bar{\eta} \in (0, \sqrt{a_\infty}) \quad\mbox{ such that
}\quad\lim_{|x|\to \infty} 
[a(x)-a_\infty] e^{\bar{\eta} 
 |x|} =  \infty, 
\eeq
used in \cite{CMcomp,CMP,CPS1,CPS2}, would not be sufficient.
}
\end{rem}

\no Next remark roughly describes the properties on which we base the
idea of the proof of Theorem \ref{T} and Proposition \ref{P1.2}, that
will be developed in Sections 4 and 5.

\begin{rem}
\label{R2.3}
{\em
The proof of Proposition \ref{P2.2} suggests that the interaction
between the points $x_1,\ldots,x_k$ tends to be attractive as
$|x_1|,\ldots,|x_k|$ tend to infinity, in the following sense.
}
\end{rem}
\no Let
$x_{1,n},\ldots,x_{k,n}$ in $\R^N$ and $\rho_n>0$
be such that $\lim_{n\to\infty}\rho_n=\infty$,
\beq
\left[{1\over k}\sum_{i=1}^k |x_{i,n}|^2\right]^{1/2}=\rho_n,\qquad
    (1+2\sigma)^{-1}\rho_n\le
|x_{i,n}|\le (1+2\sigma)\rho_n\ \mbox{ for }i=1,\ldots,k,
\eeq
and
\beq
\label{1750}
\begin{array}{rcl}
f_k(x_{1,n},\ldots,x_{k,n})&=&g^{k,\sigma}\left(\rho_n,{x_{1,n}\over
  |x_{1,n} |},\ldots,{x_{k,n}\over |x_{k,n} |}\right)\\[2ex]
&=& \max\big\{g^{k,\sigma}(\rho_n,\theta_1,\ldots,\theta_k)\
:\ (\rho_n,\theta_1,\ldots,\theta_k)\in
D^{k,\sigma}\big\}\quad\forall n\in\N.
\end{array}
\eeq
Then, for all $k\ge 2$, we have
\beq
\lim_{n\to\infty} \min\big\{|x_{i,n}-x_{j,n}|\ :\
i,j\in\{1,\ldots,k\},\ i\neq j\big\}=\infty.
\eeq
In fact, 
\beq
\liminf_{n\to\infty} \min\big\{|x_{i,n}-x_{j,n}|\ :\
i,j\in\{1,\ldots,k\},\ i\neq j\big\}<\infty
\eeq
would imply
\beq
\liminf_{n\to\infty} f_k(x_{1,n},\ldots,x_{k,n})<k  \, E_\infty(w),
\eeq
in contradiction with the fact that 
\beq
\liminf_{n\to\infty}\max \{g^{k,\sigma}(\rho_n,\theta_1,\ldots,\theta_k)\
:\ (\rho_n,\theta_1,\ldots,\theta_k)\in
D^{k,\sigma}\big\}\ge k \, E_\infty(w),
\eeq
that follows arguing exactly as in the proof of Proposition \ref{P2.2}
(see (\ref{e3.7}) and (\ref{kEinfty})).

\no Taking into account Proposition \ref{P2.2} and the definition of
$D_k$ and $D^{k,\sigma}$, we infer that $r_k$ and $\min\big\{|x^k_i|\
:\ i\in\{1,\ldots,k\}\big\}$ must tend to infinity as $k\to\infty$,
that is the interaction between the points $x^k_{1},\ldots,x^k_{k}$ tends to
be attractive.
As a consequence, because of the second equality in (\ref{2.21}), the
distances between the points ${x^k_1\over | x^k_1|},\ldots
,{x^k_k\over | x^k_k|}$ tend to be as large as possible, so these
points tend to be distributed in all of the sphere $S(0,1)$.

\no On the contrary, because of the first equality in (\ref{2.21}), the
distances between the numbers ${|x^k_1|\over r_k},\ldots
,{|x^k_k|\over r_k}$ tend to be as small as possible, so these numbers
tend to be all close to 1.

\no Taking into account the assumptions (\ref{1}) and (\ref{2}), we infer
that the number $r_k$ must be large enough so that the distances
between the points $x^k_1,\ldots,x^k_k$ tend to infinity, but not too
large, otherwise we would have
$f_k(x^k_1,\ldots,x^k_k)<\max_{D^{k,\sigma}}g^{k,\sigma}$ in
contradiction with (\ref{2.21}).


\sezione{Asymptotic estimates}


\no In this section we describe the asymptotic behaviour as
$k\to\infty$ of the sequence $((x_1^k,\ldots, x_k^k ))_k$ given by
Proposition \ref{P2.2} and of the functions $u_k$, minimizing the
energy functional $E$ in the set $S_{x_1^k,\ldots, x_k^k  }$.

\begin{prop}
\label{ktoinfty}
Assume that the potential $a(x)$ satisfies conditions (\ref{1}) and
(\ref{2}).
Let $((x_1^k,\ldots, x_k^k))_k$ be a sequence
provided by Proposition \ref{P2.2} and, for all $k\in\N$, let $u_k$ be a
function in $S_{x_1^k,\ldots, x_k^k}$ such that $E(u_k)=f_k(
x_1^k,\ldots, x_k^k )$. 

\no Then, the following properties hold:
\beq
\label{a3}
\lim_{k\to\infty} \min\{|x_i^k|\ :\ i=1,\ldots,k\}=\infty
\eeq
\beq
\label{b3}
\lim_{k\to\infty} \min\{|x_i^k-x_j^k|\ :\ i\neq j,\
i,j=1,\ldots,k\}=\infty
\eeq
\beq
\label{c3}
\lim_{k\to\infty}\sup\{|u_k(x+x_i^k)-w(x)|\ :\ |x|\le R,\
i=1,\ldots,k\}=0\quad\forall R>0.
\eeq
Moreover, there exists $\bar k>0$ such that, for all $k\ge\bar k$,
\beq
\label{eq}
-\Delta u_k(x)+a(x)u_k(x)=u^p_k(x)+\sum_{i=1}^k
(u_k)_\d^i(x)[\lambda_i^k\cdot (x-x_i^k)]\quad\forall x\in\R^N
\eeq
where $\lambda_1^k,\ldots,\lambda_k^k$ are the Lagrange
multipliers of $u_k$, and $u_k\to 0$ uniformly on the
compact subsets of $\R^N$, as $k\to\infty$.  
\end{prop}

\no\proof Property (\ref{a3}) is a direct consequence of the
definitions of $D_k$ and $D^{k,\sigma}$.

\no In order to prove (\ref{b3}) we argue by contradiction and assume
that
\beq
\liminf_{k\to\infty} \min\{|x_i^k-x_j^k|\ :\ i\neq j,\
i,j=1,\ldots,k\}<\infty.
\eeq
Without any loss of generality, we can assume that 
\beq
 \min\{|x_i^k-x_j^k|\ :\ i\neq j,\
 i,j=1,\ldots,k\}=|x^k_1-x^k_2|\qquad\forall k\ge 2.
\eeq
So (up to a subsequence) we have
$\lim_{k\to\infty}|x^k_1-x^k_2|<\infty$. 

\no Let us set
\beq
\cL_k=\max_{\theta\in S(0,1)} \min\left\{
\left|\theta-{x_i^k\over |x_i^k|}\right|\ :\ i=3,\ldots,k\right\}\qquad
\forall k\ge 3.
\eeq
We say that  $\lim_{k\to\infty}r_k\cdot\cL_k=\infty$.
In fact, arguing by contradiction, assume that (up to a subsequence)
\beq
\label{<2}
\lim\limits_{k\to\infty}r_k\cdot\cL_k<\infty.
\eeq
In this case, if we set
\beq
\d^k_i=\min\left\{\left|{x_i^k\over |x_i^k|}-{x_j^k\over |x_j^k|}\right|\
:\ j\in\{1,\ldots,k\},\ j\neq i\right\}\qquad\forall
i\in\{1,\ldots,k\},\ \forall k\ge 2,
\eeq
we must have also 
\beq
\label{A}
\lim_{k\to\infty}r_k\max\{\d^k_i\ :\ i=1,\ldots,k\}<\infty
\eeq
otherwise, for all $k\ge 3$ we could choose $\bar\theta_k\in S(0,1)$
such that  (up to a subsequence)
\beq
\lim_{k\to\infty} r_k\cdot \min\left\{
\left|\bar \theta_k-{x_i^k\over |x_i^k|}\right|\ :\
i=3,\ldots,k\right\}=\infty
\eeq
in contradiction with (\ref{<2}).

\no As a consequence of (\ref{A}), by using (\ref{a3}) and $a(x)\to
a_\infty$, and arguing as in Proposition 5.3  in \cite{CMP}, we obtain
\beq
\label{e4.12}
\liminf_{k\to\infty} {1\over k}g^{k,\sigma}\left(r_k,
{x_1^k\over |x_1^k|},\ldots,{x_k^k\over |x_k^k|}\right)
\le
\liminf_{k\to\infty} {1\over k}f_k\left(r_k\cdot
{x_1^k\over |x_1^k|},\ldots,r_k\cdot{x_k^k\over |x_k^k|}\right)
< E_\infty(w)
\eeq
in contradiction with Proposition \ref{P2.2}, which implies
\beq
\liminf_{k\to\infty} {1\over k}g^{k,\sigma}\left(r_k,
{x_1^k\over |x_1^k|},\ldots,{x_k^k\over |x_k^k|}\right)
\ge E_\infty(w).
\eeq
Thus, we have proved that $\lim_{k\to\infty}r_k\cdot
\cL_k=\infty$.
As a consequence, for all $k\ge 3$ we can choose
$\theta^k_1,\theta^k_2$ in $S(0,1)$ such that
\beq
\label{B2}
\begin{array}{c}
\lim\limits_{k\to\infty}r_k\, |\theta^k_1-\theta^k_2|=\infty \quad\mbox{ and
}
\\[2ex] 
\lim\limits_{k\to\infty}r_k\,
\min\left\{\left|\theta^k_i-{x^k_j\over|x^k_j|}\right|\ :\
  j=3,\ldots,k\right\}=\infty\quad\mbox{ for }i=1,2.
\end{array}
\eeq
Then, consider $(y^k_1,\ldots,y^k_k)$ in
$D^{k,\sigma}\left(r_k,\theta^k_1,\theta^k_2,
{x^k_3\over|x^k_3|},\ldots,{x^k_k\over|x^k_k|}\right)$ such that
\beq
f_k(y^k_1,\ldots,y^k_k)=g^{k,\sigma}\left(r_k,\theta^k_1,\theta^k_2,
{x^k_3\over|x^k_3|},\ldots,{x^k_k\over|x^k_k|}\right)
\eeq
and two points $z^k_1,z^k_2$ in $D_k$ such that
\beq
{z^k_i\over|z^k_i |}= {x^k_i\over|x^k_i |}\ \mbox{ for
}i=1,2
\eeq
and
\beq
|z_1^k|^2+|z_2^k|^2=|y_1^k|^2+|y_2^k|^2,\qquad
|\, |z_1^k|-|z_2^k|\,|\le 4\, R_\d.
\eeq
Notice that 
\beq
\label{1644}
\limsup_{k\to\infty}|z^k_1-z^k_2|<\infty
\eeq
 because
$\lim_{k\to\infty}|x^k_1-x^k_2|<\infty$.
Moreover, from (\ref{B2}) we obtain
\beq
\label{1645}
\lim_{k\to\infty}|y^k_1-y^k_2|=\infty\quad\mbox{ and
}\lim_{k\to\infty}\min\{|y^k_i-y^k_j|\ :\
j=3,\ldots,k\}=\infty\quad\mbox{ for }i=1,2.
\eeq
From (\ref{1644}) and (\ref{1645}), by the arguments used for
(\ref{e4.12}), we infer that 
\beq
\label{418}
\liminf_{k\to\infty}[f_k(y^k_1,\ldots,y^k_k)-f_k(z^k_1,z^k_2,y^k_3,\ldots,y^k_k)]>0,
\eeq
which implies
\beq
\liminf_{k\to\infty}\left[g^{k,\sigma}\left(r_k,\theta^k_1,\theta^k_2,{x^k_3\over
    |x^k_3|},\ldots,{x^k_k\over|x^k_k|}\right)
-
g^{k,\sigma}\left(r_k,{x^k_1\over
    |x^k_1|},\ldots,{x^k_k\over|x^k_k|}\right)\right]>0
\eeq
because, by the definition of $g^{k,\sigma}$ in (\ref{g}),  
\beq
f_k(z^k_1,z_2^k,y_3^k,\ldots,y_k^k)\ge g^{k,\sigma}\left(r_k,{x_1^k\over|x_1^k|},\ldots,
{x_k^k\over|x_k^k|}\right)\qquad\forall k\ge 3.
\eeq
Therefore, we have 
\beq
g^{k,\sigma}\left( r_k,\theta_1^k,\theta_2^k,{x_3^k\over|x_3^k|},\ldots,
{x_k^k\over|x_k^k|}\right)>g^{k,\sigma}\left( r_k,{x_1^k\over|x_1^k|},\ldots,
{x_k^k\over|x_k^k|}\right)
\eeq
for $k$ large enough, in contradiction with the fact that
$g^{k,\sigma}\left( r_k,{x_1^k\over|x_1^k|},\ldots,
  {x_k^k\over|x_k^k|}\right)=\max\limits_{D^{k,\sigma}}g^{k,\sigma}$. 
Thus, (\ref{b3}) is proved.

\no Property (\ref{c3}) follows from (\ref{w}) taking into account
(\ref{a3}) and (\ref{b3}).
From (\ref{c3}), (\ref{e}) and (\ref{lambda}) we deduce that, for $k$
large enough, the function $u_k$ solves the equation (\ref{eq}).
Finally, taking into account (\ref{B}) and (\ref{a3}), since $w(x)\to
0$ as $|x|\to\infty$,  from
(\ref{c3}) we deduce that
$u_k\to 0$ as $k\to\infty$ and the convergence is uniform on the
compact subsets of $\R^N$.

\qed

Let us set 
\beq
\label{Tk}
\Gamma_k=\min\left\{\left|{x_i^k\over |x_i^k|}-{x_j^k\over |x_j^k|}\right|\ :
\ i,j\in\{1,\ldots,k\},\ i\neq j\right\}\qquad\forall k\ge 2
\eeq
and
\beq
\label{Lk}
\Lambda_k=\max_{\theta\in S(0,1)}\min\left\{\left|\theta-{x_i^k\over
      |x_i^k|}\right|\ :\ i=1,\ldots,k\right\}\qquad\forall k\ge 1.
\eeq
Then, the following lemma holds.

\begin{lemma}
\label{L3.2a}
Assume that the potential $a(x)$ satisfies the conditions (\ref{1})
and (\ref{2}).
Let $(x^k_1,\ldots,x^k_k)$ and $r_k$ be as in Proposition \ref{P2.2}.
Let $\Gamma_k$ and $\Lambda_k$ be the positive numbers defined in
(\ref{Tk}) and (\ref{Lk}).
Then,
\beq
\label{Txr}
\lim_{k\to\infty}\Gamma_k=0\quad\mbox{ and }\
\lim_{k\to\infty}\Gamma_k\cdot r_k=\infty.
\eeq
If we assume in addition that condition (\ref{1.3bis}) holds, then
there exists $\tilde\sigma>0$ such that for all $\sigma\in
]0,\tilde\sigma[$ we have
\beq
\label{L/T}
\limsup_{k\to\infty}{\Lambda_k\over\Gamma_k}<\infty
\eeq
(notice that, as $x_1^k,\ldots,x_k^k$, also $\Lambda_k$ and $\Gamma_k$
depend on the parameter $\sigma$ introduced to define $D^{k,\sigma}$).
\end{lemma}

\no\proof
Notice that the balls $\left({x_1^k\over |x_1^k|},{\Gamma_k\over 2}\right),$\ldots,
$\left({x_k^k\over |x_k^k|},{\Gamma_k\over 2}\right)$ are pairwise
disjoint, so we must have $\lim_{k\to\infty}\Gamma_k=0$ because
$S(0,1)$ is a bounded set.

\no In order to prove that $\lim_{k\to\infty}\Gamma_k\cdot r_k=\infty$, we
argue by contradiction and assume that (up to a subsequence)
$\lim_{k\to\infty}\Gamma_k\cdot r_k<\infty$.
Without any loss of generality, in the following we assume also that
$\Gamma_k=\left|{x^k_1\over |x^k_1 |}- {x^k_2\over |x^k_2 |}\right|$
    $\forall k\ge 2$.

\no Then, consider two points $z_1^k$ and $z_2^k$ in $D_k$ such that 
\beq
{z^k_i\over |z^k_i|}={x^k_i\over |x_i^k|}\ \mbox{ for }i=1,2
\eeq
and
\beq
|z^k_1|^2+|z^k_2|^2=|x^k_1|^2+|x_2^k|^2,\qquad |\,
|z_1^k|-|z^k_2|\,|<4\, R_\d.
\eeq
Notice that $\lim_{k\to\infty} \Gamma_k\cdot r_k<\infty$ implies
  $\limsup_{k\to\infty}|z^k_1-z^k_2|<\infty$.
Therefore, taking into account (\ref{a3}) and (\ref{b3}) and arguing
as for (\ref{418}), we obtain
\beq
\liminf_{k\to\infty}[f_k(x_1^k,\ldots,x_k^k)-f_k(z_1^k,z_2^k,x_3^k,\ldots,x_k^k)]>0
\eeq
which is a contradiction because
$f_k(x_1^k,\ldots,x_k^k)=g^{k,\sigma}\left(r_k,{x_1^k\over
    |x_1^k|},\ldots,{x^k_k\over |x_k^k|}\right)$.
Thus (\ref{Txr}) is proved.

\no For the proof of (\ref{L/T}) we argue again by contradiction and
assume that (up to a subsequence)
\beq
\label{infty3}
\lim_{k\to\infty}{\Lambda_k\over\Gamma_k}=\infty.
\eeq
Therefore, we can choose $\theta_1^k$, $\theta_2^k$ in $S(0,1)$ such
that
\beq
\label{theta}
\begin{array}{c}
\lim\limits_{k\to\infty}{|\theta_1^k-\theta_2^k|\over
  \Gamma_k}=\infty\quad\mbox{ and }
\\[2ex]
 \lim\limits_{k\to\infty}{1\over
  |\theta_1^k-\theta_2^k|}\cdot\min\left\{\left|\theta_i^k-{x_j^k\over
      |x^k_j|}\right|\ :\ j=3,\ldots,k\right\}=\infty\ \mbox{ for }\
i=1,2.
\end{array}
\eeq
Now, consider $(y^k_1,\ldots,y^k_k)$ in
$D^{k,\sigma}\left(r_k,\theta^k_1,\theta_2^k,{x^k_3\over |x^k_3|},\ldots,
{x^k_k\over |x^k_k|}\right)$ such that 
\beq
\label{f=g}
f_k(y^k_1,\ldots,y^k_k)=g^{k,\sigma}\left(r_k,\theta^k_1,\theta^k_2,{x^k_3\over |x^k_3|},\ldots,
{x^k_k\over |x^k_k|}\right)
\eeq
and the points $\xi^k_1$, $\xi_2^k$ in $D_k$ such that
\beq
\label{4.34}
|\xi^k_i|=|y^k_i|,\qquad {\xi^k_i\over |\xi^k_i|}={x^k_i\over
  |x^k_i|}\qquad\mbox{ for }\ i=1,2.
\eeq
We claim that 
\beq
\label{>xi}
f_k(y^k_1,\ldots,y^k_k)>f_k(\xi^k_1,\xi^k_2,y^k_3,\ldots,y^k_k)
\eeq
for $k$ large enough.
Once (\ref{>xi}) is proved, we are done. 
Indeed, it  is a contradiction   because 
\beq
f_k(\xi_1^k,\xi_2^k,y^k_3,\ldots,y^k_k)\ge
g^{k,\sigma}\left(r_k,{x^k_1\over|x^k_1|},\ldots,{x^k_k\over
    |x^k_k|}\right)=\max_{D^{k,\sigma}} g^{k,\sigma}
\eeq
and (\ref{f=g}) holds.

\no In order to prove (\ref{>xi}), let us set
\beq
\e(x_1,x_2)=f_1(x_1)+f_1(x_2)-f_2(x_1,x_2)\qquad\forall (x_1,x_2)\in
D_2
\eeq
and
\beq
\e(x_1,\ldots,x_k)=f_2(x_1,x_2)+f_{k-2}(x_3,\ldots,x_k)-f_k(x_1,\ldots,x_k)\qquad\forall
(x_1,\ldots,x_k)\in D_k.
\eeq
Then we obtain
\beq
\label{abeta}
\begin{array}{rcl}
\vspace{2mm}
f_k(\xi^k_1,\xi^k_2,y^k_3,\ldots,y^k_k) &=&
f_k(y^k_1,y_2^k,y_3^k,\ldots,y^k_k)
-\e(\xi^k_1,\xi^k_2)+\e(y^k_1,y_2^k )\\
\vspace{2mm}
& & 
+\e(y^k_1,y_2^k,y_3^k,\ldots,y^k_k)-\e(\xi^k_1,\xi^k_2,y^k_3,\ldots,y^k_k)
\\
& & +f_1(\xi_1^k)+f_1(\xi_2^k)-f_1(y_1^k)-f_1(y_2^k).
\end{array}
\eeq
Arguing as in the proof of Proposition 4.5 in \cite{CMP}, one can
verify that $\e(\xi^k_1,\xi^k_2)$, $\e(y^k_1,y_2^k )$,
$\e(\xi^k_1,\xi^k_2,y^k_3,\ldots,y^k_k)$ and
$\e(y^k_1,y_2^k,y_3^k,\ldots,y^k_k)$ are positive numbers, for
all $k>2$, that tend to zero as $k\to\infty$ and there exists a
suitable positive constant $\bar c$ such that 
\beq
\label{alpha}
\liminf_{k\to\infty}\e(\xi^k_1,\xi^k_2)\, e^{\bar c\, |\xi^k_1-\xi^k_2|}>0.
\eeq
Moreover, from (\ref{4.34}) and assumption (\ref{1.3bis}) we infer that
\beq
\label{1211}
|f_1(\xi_1^k)+f_1(\xi_2^k)-f_1(y_1^k)-f_1(y_2^k)|\le \tilde c e^{-\tilde\eta\,
  {r_k\over 1+2\sigma}} 
\eeq
for a suitable constant $\tilde c>0$.

\no Now, let us prove that (\ref{theta}) implies
\beq
\label{star}
\lim_{k\to\infty}(|y^k_1-y^k_2|-|\xi^k_1-\xi^k_2|)=\infty.
\eeq
First notice that, since $|\xi_i^k|=|y^k_i|$ for $i=1,2$, we get
\beq
\label{2star}
\begin{array}{rcl}
|y^k_1-y^k_2|-|\xi^k_1-\xi^k_2|&=&{|y^k_1-y^k_2|^2-|\xi^k_1-\xi^k_2|^2\over
|y^k_1-y^k_2|+|\xi^k_1-\xi^k_2|}
\\[1.5ex]
& = &
{2\big(\xi^k_1\cdot \xi^k_2-y^k_1\cdot y^k_2\big)\over|y^k_1-y^k_2|+|\xi^k_1-\xi^k_2|}
\\[1.5ex]
& = &
{2|\xi_1^k|\,
  |\xi^k_2|\cdot\left[{x^k_1\over|x^k_1|}\cdot
    {x^k_2\over|x^k_2|}-\theta^k_1\cdot \theta^k_2\right]\over|y^k_1-y^k_2|+|\xi^k_1-\xi^k_2|}\qquad\forall
  k\ge 2.
\end{array}
\eeq
Moreover, since $\Gamma_k=\left|{x^k_1\over |x^k_1|}-{x^k_2\over
    |x^k_2|}\right|$, we can write
\beq
\label{438}
{x^k_1\over|x^k_1|}\cdot {x^k_2\over|x^k_2|}-\theta^k_1\cdot \theta^k_2
= {1\over 2}(|\theta^k_1-\theta^k_2|^2-\Gamma_k^2).
\eeq

\no Since $\theta^k_1$ and $\theta^k_2$ satisfies (\ref{theta}), we
have $\lim\limits_{k\to\infty}|\theta^k_1-\theta^k_2|=0$. 
Moreover we say that
\beq
\label{limsup3}
\limsup_{k\to\infty} {|y^k_1-y^k_2|\over r_k\,
  |\theta^k_1-\theta^k_2|}<\infty.
\eeq
In fact, arguing by contradiction, assume that (up to a subsequence)
\beq
\label{y/rtheta}
\lim_{k\to\infty} {|y^k_1-y^k_2|\over r_k\,
  |\theta^k_1-\theta^k_2|}=\infty
\eeq
and consider two points $\psi^k_1$, $\psi^k_2$ in $D_k$ such that
\beq
 {\psi^k_i\over |\psi^k_i|}={y^k_i\over
  |y^k_i|}=\theta^k_i\qquad\mbox{ for }\ i=1,2
\eeq
and\beq
|\psi^k_1|^2+|\psi^k_2|^2=|y^k_1|^2+|y^k_2|^2,\qquad
|\,|\psi^k_1|-|\psi^k_2|\, |\, \le 4R_\d.
\eeq
Then, we can argue as in the proof of (\ref{418}) and we infer  from
(\ref{theta}) and (\ref{y/rtheta}) that, for $k$ large enough, 
\beq
f_k(y^k_1,\ldots,y^k_k)>f_k(\psi^k_1,\psi^k_2,y^k_3,\ldots,y^k_k)
\eeq
in contradiction with (\ref{f=g}).
Thus (\ref{limsup3}) holds and, as a consequence, we have also
\beq
\label{xi/rtheta}
\limsup_{k\to\infty} {|\xi^k_1-\xi^k_2|\over r_k\,
  |\theta^k_1-\theta^k_2|}<\infty.
\eeq
Hence, since  (\ref{2star}), (\ref{438}) and (\ref{theta}) give
\beq
\liminf_{k\to\infty} (|y^k_1-y^k_2|-|\xi^k_1-\xi^k_2|)\ge
\liminf_{k\to\infty}{(1+2\sigma)^{-2}\cdot
  r^2_k|\theta^k_1-\theta^k_2|^2\over
  |y^k_1-y^k_2|+|\xi^k_1-\xi^k_2|},
\eeq
we obtain (\ref{star}) by  (\ref{limsup3}), (\ref{xi/rtheta}) and
\beq
\lim_{k\to\infty}r_k|\theta^k_1-\theta^k_2|\ge
\lim_{k\to\infty}r_k\Gamma_k=\infty.
\eeq
From (\ref{star}) it follows that 
\beq
\label{betak}
\lim_{k\to\infty}{\e(y_1^k,y_2^k)\over\e(\xi^k_1,\xi^k_2)}=0.
\eeq
In analogous way, from the choice of $\theta_1^k,\theta_2^k$ it follows that 
\beq
\label{n*}
\lim_{k\to\infty} {\e(\xi^k_1,\xi^k_2,y_3^k,\ldots,y_k^k)\over
  \e(\xi^k_1,\xi^k_2)}=
\lim_{k\to\infty} {\e(y_1^k,y_2^k,\ldots,y_k^k)\over
  \e(\xi^k_1,\xi^k_2)}=0.
\eeq
Let us recall that $\xi^k_1,\xi^k_2$ depend on the choice of $\sigma$.
Now, we prove that there exists $\tilde \sigma>0$ such that, for all
$\sigma\in]0,\tilde\sigma[$,  
\beq
\label{tilde}
\lim_{k\to\infty}e^{-\tilde\eta\, {r_k\over 1+2\sigma}}e^{\bar c\,
  |\xi^k_1-\xi^k_2|}=0.
\eeq
Let us choose $\tilde \sigma>0$ small enough so that $\bar
c[(1+2\tilde\sigma)^2-1]<\tilde\eta$.
Then, since $\Gamma_k\to 0$, for all $\sigma\in]0,\tilde\sigma[$ we
have
\beq
\limsup_{k\to\infty}\bar c\, {1+2\sigma \over r_k} \cdot |\xi^k_1-\xi^k_2|\le \bar
c[(1+2\sigma)^2-1]\le\bar c[(1+2\tilde
\sigma)^2-1]<\tilde\eta
\eeq
which implies (\ref{tilde}).
Hence, from (\ref{alpha}), (\ref{betak}), (\ref{n*}), (\ref{1211}), (\ref{tilde}) and
(\ref{abeta}) we obtain our claim (\ref{>xi}), so the proof is complete.

\qed

\no As it is specified in next corollary, property
(\ref{L/T}) implies that the points ${x^k_1\over
  |x^k_1|},\ldots,{x^l_k\over |x_k^k|}$ tend, as $k\to\infty$, to
spread on all of $S(0,1)$ and that the limit density of distribution
is everywhere positive on $S(0,1)$.

\begin{cor}
\label{Cnew}
Assume that all the conditions of Lemma \ref{L3.2a} are satisfied.
Then $\lim\limits_{k\to\infty}\Lambda_k$ $=0$.

\no Moreover, if for all $x\in S(0,1)$ and $r>0$ we denote by
$N_k(x,r)$ the number of elements of the set $\Bigl\{x^k_i$ :
$i\in\{1,\ldots,k\}$, ${x^k_i\over |x^k_i|}\in B(x,r)\Bigr\}$, then there
exists $c>0$ such that
\beq
\label{Nk}
\liminf_{\e\to 0}\liminf_{k\to\infty} {N_k(x,\e)\over k\cdot \e^{N-1}}\ge
c\qquad\forall x\in S(0,1).
\eeq
\end{cor}

\no \proof
Since $\lim_{k\to\infty}\Gamma_k=0$, $\lim_{k\to\infty}\Lambda_k=0$
follows directly from (\ref{L/T}).
In order to prove (\ref{Nk}) we argue by contradiction and assume that
there exists a sequence $(x_n)_n$ in $S(0,1)$ such that
\beq
\label{0}
\lim_{n\to\infty}\liminf_{\e\to 0}\liminf_{k\to\infty} 
{N_{k}(x_n,\e)\over k\cdot \e^{N-1}}
=0.
\eeq
Taking into account the definition of $\Lambda_k$ and $\Gamma_k$, we
infer that
\beq
\inf\left\{ N_{k}(x_n,\e)\cdot\left({\Lambda_{k}\over
    \e}\right)^{N-1}\ :\  \e>0,\ \ k, n\in\N\right\}>0
\eeq
and
\beq
\sup\left\{k\, \Gamma_k^{N-1}\ :\ k\in\N\right\}<\infty.
\eeq
It follows that, for a suitable positive constant $c$,
\beq
\begin{array}{rl}
\vspace{2mm}
\limsup\limits_{k\to\infty}{\Lambda_{k}\over \Gamma_{k}}
&\ge
c\,\lim\limits_{n\to\infty}\limsup\limits_{\e\to 0}\limsup\limits_{k\to\infty} 
 \left[{k\,\e^{N-1}\over
    N_{k}(x_n,\e)}\right]^{1\over N-1}\\
&=c\,\lim\limits_{n\to\infty}\left[\liminf\limits_{\e\to 0}\liminf\limits_{k\to\infty} 
 {{N_{k}(x_n,\e)}\over k\,\e^{N-1}}
    \right]^{-{1\over N-1}}=\infty
\end{array}
\eeq
in contradiction with (\ref{L/T}).
Thus (\ref{Nk}) is proved.

\qed

\begin{lemma}
\label{sigma}
Assume that the potential $a(x)$ satisfies conditions (\ref{1}), (\ref{2}) and
(\ref{1.3bis}).
For all $k\ge 2$, let $r_k$ and $(x_1^k,\ldots,x_k^k)$ satisfy all the properties
described  in Proposition \ref{P2.2}.
Let $\tilde\sigma$ be as in Lemma \ref{L3.2a}. 
Then, there exists $\bar\sigma\in]0,\tilde\sigma[$ such that, for $\sigma=\bar\sigma$,
\beq
\label{e3.5}
\lim_{k\to\infty}{1\over r_k}\min\{|x_i^k|\ :\
i=1,\ldots,k\}
=
\lim_{k\to\infty}{1\over r_k}\max\{|x_i^k|\ :\
i=1,\ldots,k\}= 1.
\eeq
\end{lemma}

\no\proof
Let us consider a sequence  $(\bar\sigma_n)_n$ in $]0,\tilde\sigma[$ such
that $\lim_{n\to\infty}\bar\sigma_n=0$.
We shall prove that there exists $\bar n\in\N$ such that the assertion
of the lemma holds for $\bar\sigma=\bar\sigma_{\bar n}$.
Let us recall that $(x_1^k,\ldots,x_k^k)$ and $r_k$ depend also on
$\sigma$.
Therefore, if $\sigma=\bar\sigma_n$, in this proof we write, more
explicitly, $(x_1^{k,\bar\sigma_n},\ldots,x_k^{k,\bar\sigma_n})$ and
$r_{k,\bar\sigma_n}$ $\forall n\in\N$.

\no Since 
\beq
\label{3.6}
\min\{|x_i^{k,\bar\sigma_n}|\ :\ i=1,\ldots,k\}\le
r_{k,\bar\sigma_n}\le
\max\{|x_i^{k,\bar\sigma_n}|\ :\ i=1,\ldots,k\},\quad\forall k\in\N,\
\forall n\in\N,
\eeq
we have to prove that 
\beq
\limsup_{k\to\infty}{1\over
  r_{k,\bar\sigma_n}}\max\{|x_i^{k,\bar\sigma_n}|\ :\
i=1,\ldots,k\}\le 1
\le
\liminf_{k\to\infty}{1\over
  r_{k,\bar\sigma_n}}\min\{|x_i^{k,\bar\sigma_n}|\ :\
i=1,\ldots,k\}
\eeq
for some $n\in\N$.
Arguing by contradiction, assume that, for all $n\in\N$,
\beq
\label{a6}
\liminf_{k\to\infty}{1\over r_{k,\bar\sigma_n}}\min\{|x_i^{k,\bar\sigma_n}|\ :\
i=1,\ldots,k\} <1
\eeq
or
\beq
\label{b6}
\limsup_{k\to\infty}{1\over r_{k,\bar\sigma_n}}\max\{|x_i^{k,\bar\sigma_n}|\ :\
i=1,\ldots,k\}> 1.
\eeq
Let us consider, for example, the case that (\ref{b6}) is true
(similar arguments work when (\ref{a6}) holds).

\no Then,  consider the sequence of positive numbers $(\sigma_n)_n$
such that
\beq
\label{3.10}
1+2\sigma_n=\limsup_{k\to\infty}{1\over r_{k,\bar\sigma_n}}\max\{|x_i^{k,\bar\sigma_n}|\ :\
i=1,\ldots,k\} \qquad\forall n\in\N.
\eeq
Notice that $\lim_{n\to\infty}\sigma_n=0$ because
$\sigma_n\le\bar\sigma_n$ $\forall n\in\N$.
Now, for all $n\in\N$, consider the set
\beq
\label{Vn}
V_n=\{\tau\in \Z^N\ :\ (\tau+[-1,1]^N)\cap
S(0,1/\sigma_n)\neq\emptyset\}
\eeq
 and denote by $\nu_n$ the number of elements of $V_n$.
It is clear  that $S(0,1/\sigma_n)\subseteq\cup_{\tau\in
  V_n}(\tau+[-1,1]^N)$, that $\nu_n<\infty$ $\forall n\in\N$ and that
$\lim_{n\to\infty}\nu_n=\infty$.

\no Now, consider a sequence $(\gamma_n)_n$ such that $\gamma_n>0$
$\forall n\in\N$ and
\beq
\label{gamma}
\lim_{n\to\infty}{\nu_n\over\gamma_n}=0.
\eeq
From (\ref{3.10}) it follows that there exists a sequence $(k_n)_n$ in
$\N$ such that 
\beq
\label{b7}
{1\over r_{k_n,\bar\sigma_n}}\max\{|x_i^{k_n,\bar\sigma_n}|\ :\
i=1,\ldots,k_n\}>{ 1+\sigma_n},\quad\forall n\in\N
\eeq
and in addition
\beq
\label{kn}
\lim_{n\to\infty}{1\over k_n}\,
\left({\gamma_n\over\sigma_n}\right)^{N-1}=0.
\eeq

\no Taking into account Corollary \ref{Cnew}, up to a subsequence we have
\beq
\label{e3.15}
\lim_{\e\to 0}\lim_{n\to\infty}{N_{k_n}(x,\e)\over
  k_n\e^{N-1}}\ge c>0\qquad\forall x\in S(0,1)
\eeq
(where $c$ is a positive constant independent of $x$).
Moreover, taking also into account (\ref{kn}), we infer that
\beq
\liminf_{n\to\infty}\inf\left\{N_{k_n}\left(x,r\,
    {\sigma_n\over\gamma_n}\right)\ :\ x\in
  S(0,1)\right\}= \infty\qquad\forall r>0.
\eeq

\no Now, let us set 
\begin{eqnarray}\vspace{2mm}
M_n={\gamma_n\over \sigma_n r_{k_n,\bar\sigma_n}}\, \max_{\tau\in V_n}\max\bigg\{
|x_i^{k_n,\bar\sigma_n}|-|x_j^{k_n,\bar\sigma_n}| & : & i,j=1,\ldots,k_n, \\
\nonumber & &\hspace{-2cm}\left. {x_i^{k_n,\bar\sigma_n}\over
  |x_i^{k_n,\bar\sigma_n} |} \mbox{ and } {x_j^{k_n,\bar\sigma_n}\over
  |x_j^{k_n,\bar\sigma_n} |}\ \mbox{ in }\sigma_n(\tau+[-1,1]^N)\right\}.
\end{eqnarray}
From the definition of the function $g^{k,\sigma}$ and the properties
of the points
$x_1^{k_n,\bar\sigma_n}, \ldots,$ $x_{k_n}^{k_n,\bar\sigma_n}$
described in Remark \ref{R2.3} and in Corollary \ref{Cnew}, since  
the curvature of the sphere $S(0,1/\sigma_n)$ tends to zero as
$n\to\infty$, we infer that $\lim_{n\to\infty}M_n=0$
otherwise, arguing as in the proofs of Proposition \ref{ktoinfty} and Lemma
\ref{L3.2a},  we would obtain
\beq
f_{k_n}(x^{k_n,\bar\sigma_n}_1,\ldots,x^{k_n,\bar\sigma_n}_{k_n})
>f_{k_n}\left(r_{k_n}\, 
  {x^{k_n,\bar\sigma_n}_1\over |x^{k_n,\bar\sigma_n}_1 |},\ldots,r_{k_n}\,
  {x^{k_n,\bar\sigma_n}_{k_n}\over
    |x^{k_n,\bar\sigma_n}_{k_n}|}\right)
\eeq
for $n$ large enough, in contradiction with the fact that
\beq
f_{k_n}(x^{k_n,\bar\sigma_n}_1,\ldots,x^{k_n,\bar\sigma_n}_{k_n})=g^{k_n,\bar\sigma_n}\left(r_{k_n},
{x^{k_n,\bar\sigma_n}_1\over |x^{k_n,\bar\sigma_n}_1 |},\ldots,
  {x^{k_n,\bar\sigma_n}_{k_n}\over
    |x^{k_n,\bar\sigma_n}_{k_n}|}\right).
\eeq
It follows that
\beq
\begin{array}{l}\vspace{2mm}
\hspace{-4mm}\limsup\limits_{n\to\infty}{1\over \sigma_n
  r_{k_n,\bar\sigma_n}}\max\{||x_i^{k_n,\bar\sigma_n}|-r_{k_n,\bar\sigma_n}|\ :\ i=1,\ldots,k_n\} \le
\\ \vspace{2mm}
\limsup\limits_{n\to\infty}{1\over \sigma_n
  r_{k_n,\bar\sigma_n}}\max\{|x_i^{k_n,\bar\sigma_n}|-|x_j^{k_n,\bar\sigma_n}|\ :\
i,j=1,\ldots,k_n\} \le
\lim\limits_{n\to\infty}{\nu_n M_n\over\gamma_n} =  0
\end{array}
\eeq
(where the last inequality holds because $S(0,1/\sigma_n)$ is a
connected set).

\no It is clear that we have a contradiction because from (\ref{b7})
we obtain 
\beq
{1\over \sigma_n r_{k_n,\bar\sigma_n}}\,\max\{|x_i^{k_n,\bar\sigma_n}|-r_{k_n,\bar\sigma_n}\ :\
i=1,\ldots,k_n\}>1\quad\forall n\in\N,
\eeq
which implies
\beq
\liminf_{n\to\infty}{1\over \sigma_n
  r_{k_n,\bar\sigma_n}}\max\{||x_i^{k_n,\bar\sigma_n}|-r_{k_n,\bar\sigma_n}|\ :\ i=1,\ldots,k_n\}\ge 1.
\eeq
Thus, the proof is complete.

\qed


\sezione{Proof of the main result and final remarks}


\no Let us denote by $\lambda_1^k$, \ldots , $\lambda_k^k$ the
Lagrange multipliers corresponding to a minimizing function $u_k$ for
the energy functional $E$ in $S_{x_1^k,\ldots,x_k^k}$ (see
(\ref{lambda})).
In order to prove that $u_k$ is a positive solution of problem
(\ref{.}), it remains to show that, for $k$ large enough,
$\lambda_i^k=0$ $\forall i\in\{1,\ldots,k\}$.
\begin{prop}
\label{P4.1}
Assume that $\sigma =\bar\sigma$ (see Lemma \ref{sigma}).
Let $r_k$ and $(x^k_1,\ldots,x^k_k)$ be as in Proposition \ref{P2.2}.
Then, there exists $\bar k'\in\N$ such that, for all $k\ge \bar k'$, 
\beq
\label{lambdamu}
{\l^k_i\cdot x^k_i\over
  |x^k_i|^2}\int_{B(x^k_i,R_\d)}[u^\d_k(x)]^2dx=\mu_k\qquad\forall
  i\in\{1,\ldots,k\}
\eeq
for a suitable $\mu_k\in\R$.
\end{prop}

\no \proof 
Notice that every function $u_k\in S_{x^k_1,\ldots,x^k_k}$, such
that $E(u_k)=f_k(x^k_1,\ldots,x^k_k)$, satisfies
\beq
\label{min2}
\begin{array}{rcl}
\vspace{2mm}
E(u_k)\hspace{-2mm}&=\min\{  E(u)& :\ u^\d\in H^1_0(\cup_{i=1}^k
B(x_i,R_\d)),\ u\in S_{y_1,\ldots,y_k},
\\ \vspace{2mm} & &
\left[{1\over k}\sum_{i=1}^k|y_i|^2\right]^{1/2}=r_k,\ y_i\in B(x_i,R_\d), 
\\ \vspace{2mm} & &
(1+2\bar\sigma)^{-1}r_k\le |y_i|\le (1+2\bar\sigma)r_k,\ {y_i\over
  |y_i|}={x^k_i\over|x^k_i|}\ \mbox{ for }i=1,\ldots,k\}.
\end{array}
\eeq
Moreover, (\ref{w}) and Lemma \ref{sigma} imply that there exists
$\bar k'\in\N$ such that, for all $k\ge \bar k'$,
\beq
\label{<}
(1+2\bar\sigma)^{-1}r_k < |x_i| < (1+2\bar\sigma)r_k \qquad \forall i\in\{1,\ldots,k\},
\eeq
and, for $i=1,\ldots,k$, there exists $w_{k,i}\in H^1_0(B(x^k_i,
R_\d))$ such that $\beta_i'(u_k)[w_{k,i}]=x^k_i$. 

\no
Therefore, arguing as in  Proposition 3.5 in \cite{CMP} we infer that for all
$k\ge \bar k'$ there exists a Lagrange multiplier $\mu_k\in\R$ such 
that
\beq
\label{muk}
E'(u_k)[w_{k,i}]={\mu_k\over 2}\, x^k_i\cdot
\beta'_i(u_k)[w_{k,i}]\qquad\forall i\in\{1,\ldots,k\}.
\eeq

\no On the other hand, from (\ref{2.5bis}) and (\ref{lambda}) we
obtain
\beq
E'(u_k)[w_{k,i}]={1\over 2}\,\lambda^k_i\cdot \beta'_i(u_k)[w_{k,i}]
\int_{B(x^k_i,R_\d)}[u^\d_k(x)]^2dx
\eeq
which, combined with (\ref{muk}), implies
\beq
\left[{\lambda^k_i\over 2} \int_{B(x^k_i,R_\d)}[u^\d_k(x)]^2dx-{\mu_k\over 2}\, x^k_i\right]\cdot
  \beta'_i(u_k)[w_{k,i}]=0,
\eeq
that is (\ref{lambdamu}).

\qed

\begin{rem}
{\em
 In the proof of next propositions we obtain some integrals of the form 
\beq
\int_{\R^N}w^\d(x)\, (D w(x)\cdot\tau)\, x\, dx
\eeq
when we apply Lemma \ref{Lw}. Let us remark that 
\beq
\label{form}
\begin{array}{rcl}
\vspace{2mm}
\int_{\R^N} w^\d(x)\, (D w(x)\cdot\tau)\, x\, dx & = &
{1\over 2}\int_{\R^N} (D[w^\d(x)]^2\cdot\tau)\, x\, dx 
\\ & = & 
-{\tau\over 2}\int_{\R^N}[w^\d(x)]^2 dx,
\end{array}
\eeq
as one can verify by direct computation. 
}
\end{rem}

\begin{prop}
\label{P4.2}
Assume that $\sigma =\bar\sigma$ (see Lemma \ref{sigma}).
Let $r_k$ and $(x^k_1,\ldots,x^k_k)$ be as in Proposition \ref{P2.2}.
Then, there exists $\bar k''\in\N$ such that
\beq
\label{lambdax}
\l^k_i\cdot x^k_i=0\qquad\forall k\ge \bar k'',\ \forall i\in\{1,\ldots,k\}.
\eeq
\end{prop}

\no \proof
From Proposition \ref{P4.1} it follows that there exists $\bar
k'\in\N$ such that (\ref{lambdamu}) holds for all $k\ge \bar k'$.
Thus, it remains to show that there exists $\bar k''>\bar k'$ such
that $\mu_k=0$ $\forall k\ge \bar k''$.
Arguing by contradiction, assume that there exists a sequence
$(k_n)_n$ in $\N$ such that $\lim_{n\to\infty}k_n=\infty$ and
$\mu_{k_n}\neq 0$ $\forall n\in\N$.
Up to a subsequence, $|\mu_{k_n}|^{-1}\mu_{k_n}\to\bar\mu$ as
$n\to\infty$, for a suitable $\bar\mu\in\{-1,1\}$.
Now, choose a sequence $(\bar\e_n)_n$ of positive numbers such that
$\lim_{n\to\infty}\bar\e_n/\mu_{k_n}=0$ (notice that, as
a consequence, $\lim_{n\to\infty}\bar\e_n=0$ because
$\lim_{n\to\infty}\mu_{k_n}=0$ as follows from
(\ref{lambdan})).
Then, set $\rho_n=(r^2_{k_n}+\bar\mu\, \bar\e_n)^{1/2}$ $\forall n\in\N$ and notice
that $\left(\rho_n,{x^{k_n}_1\over |x^{k_n}_1 |},\ldots,
  {x^{k_n}_{k_n}\over  |x^{k_n}_{k_n} |}\right)\in D^{k_n,\bar\sigma}$
for $n$ large enough and
\beq
\label{le}
g^{k_n,\bar\sigma}
\left(\rho_n,{x^{k_n}_1\over |x^{k_n}_1 |},\ldots,
  {x^{k_n}_{k_n}\over  |x^{k_n}_{k_n} |}\right)\le
g^{k_n,\bar\sigma}
\left(r_{k_n},{x^{k_n}_1\over |x^{k_n}_1 |},\ldots,
  {x^{k_n}_{k_n}\over  |x^{k_n}_{k_n} |}\right).
\eeq
Let us choose $(\bar y^{k_n}_1,\ldots,\bar y^{k_n}_{k_n})$ in
$D_{k_n}$ satisfying
\beq
\left[{1\over {k_n}}\sum_{i=1}^{k_n}|\bar
  y^{k_n}_i|^2\right]^{1/2}=\rho_n,
\eeq
\beq
(1+2\bar\sigma)^{-1}\rho_n\le |\bar
y^{k_n}_i|\le(1+2\bar\sigma)\rho_n,\quad
{\bar y^{k_n}_i\over |\bar y^{k_n}_i |}=
{x^{k_n}_i\over | x^{k_n}_i |}\quad\mbox{ for
}i=1,\ldots,k_n
\eeq
and
\beq
f_{k_n}(\bar y^{k_n}_1,\ldots,\bar y^{k_n}_{k_n} )=
g^{k_n,\bar\sigma}\left(\rho_n,{ x^{k_n}_1\over | x^{k_n}_1|},
\ldots,
{ x^{k_n}_{k_n}\over | x^{k_n}_{k_n}|}\right).
\eeq
Then, we have 
\beq
f_{k_n}(\bar y^{k_n}_1,\ldots,\bar y^{k_n}_{k_n} )\le
f_{k_n}( x^{k_n}_1,\ldots, x^{k_n}_{k_n} )
\eeq
which implies $E(\bar v_n)\le E(u_{k_n})$ for every $\bar v_n\in
S_{\bar y^{k_n}_1,\ldots,\bar y^{k_n}_{k_n} }$ such that $E(\bar
v_n)=f_{k_n}(\bar y^{k_n}_1,\ldots,$ $\bar y^{k_n}_{k_n} )$.
As in the proof of Theorem 1.1 in \cite{CPS1}, we have
\beq
\label{4.15}
\liminf_{n\to\infty}{E(\bar v_n)-E(u_{k_n})\over \bar
  \e_n|\mu_{k_n}|\, k_n}\ge\liminf_{n\to\infty}{1\over \bar
  \e_n|\mu_{k_n}|\, k_n}\cdot\sum_{i=1}^{k_n}C_{n,i}
\eeq
where
\beq
C_{n,i}=\int_{B(x_i^{k_n},R_\delta)}u^\d_{k_n}(x)\, [\bar
v_n(x)-u_{k_n}(x)]\, [\lambda_i^{k_n}\cdot(x-x_i^{k_n})]\, dx.
\eeq
For every $n\in \N$, let us choose $i_n\in\{1,\ldots,k_n\}$
and assume that (up to a subsequence) ${x_{i_n}^{k_n}\over |
  x_{i_n}^{k_n}|}\to\bar x\in S(0,1)$ and ${1\over \bar\e_n}(\bar
y_{i_n}^{k_n}-x_{i_n}^{k_n})\cdot x_{i_n}^{k_n}\to\bar c\in\R$ as
$n\to\infty$.
Then, taking into account Lemma \ref{Lw}, we obtain 
\beq
\label{4.18}
 \liminf_{n\to\infty} {C_{n,i_n}\over |\mu_{k_n}|\, \bar\e_n}\ge
 -
\bar c \left[\int_{B(0,R_\d)}[w^\d (x)]^2dx\right]^{-1}
 \int_{B(0,R_\delta)}w^\d(x) \, (D
w(x)\cdot \bar x)\, (x\cdot\bar\mu\bar x)\, dx.
\eeq
Notice that in (\ref{4.18}) the integral does not depend on $\bar x$,
in the sense that its value remains unchanged if we replace $\bar x$
by any other $\bar x'\in S(0,1)$.
Therefore, it follows that
\beq
\begin{array}{l}\vspace{2mm}
\hspace{-2mm}
\liminf\limits_{n\to\infty}
{1\over \bar\e_n|\mu_{k_n}|\, k_n}
\sum\limits_{i=1}^{k_n} C_{n,i}\ge 
-\left[\int_{B(0,R_\d)}[w^\d(x)]^2dx\right]^{-1}\cdot
\\
\hspace{5mm}\cdot \int_{B(0,R_\d)}w^\d(x) \, (D
w(x)\cdot\bar x)\, (x\cdot\bar\mu\bar x)\, dx\cdot
\liminf\limits_{n\to\infty}{1\over \bar\e_n k_n}
\sum\limits_{i=1}^{k_n}(\bar y_i^{k_n}-x_i^{k_n})\cdot x_i^{k_n}.
\end{array}
\eeq
On the other hand, we have 
\beq
\lim_{n\to\infty}{1\over \bar\e_n k_n} \sum\limits_{i=1}^{k_n}(\bar
y_i^{k_n}-x_i^{k_n})\cdot x_i^{k_n}={\bar\mu\over 2}
\eeq
because of the definition of $\rho_n$.
Therefore, using also (\ref{form}) we obtain
\beq
\label{ge}
\begin{array}{l}\vspace{2mm}
\hspace{-15mm} 
\liminf_{n\to\infty}
{E(\bar v_n)-E(u_{k_n})\over \bar\e_n|\mu_{k_n} |\, k_n }\ge -
{1\over 2} \left[\int_{B(0,R_\d)}[w^\d(x)]^2dx\right]^{-1}
\cdot
\\
\hspace{15mm}\cdot 
\int_{B(0,R_\d)}w^\d(x)\, (Dw(x)\cdot\bar
x)(x\cdot\bar x)\, dx>0, 
\end{array}
\eeq
\no which gives a contradiction because $E(\bar v_n)\le E(u_{k_n})$
for $n$ large enough.   
So the proof is complete.

\qed

\begin{prop}
\label{P4.3}
Assume that $\sigma =\bar\sigma$ (see Lemma \ref{sigma}).
Let $r_k$ and $(x^k_1,\ldots,x^k_k)$ be as in Proposition \ref{P2.2}.
Then, there exists $\bar k\in\N$ such that
\beq
\label{lambda=0}
\lambda_i^k=0
\qquad\forall k\ge\bar k,\ \forall i\in\{1,\ldots,k\}.
\eeq
\end{prop}

\no\proof
Arguing by contradiction, assume that there exists a sequence
$(k_n)_n$ in $\N$ such that $\lim_{n\to \infty}k_n=\infty$ and,
for all $n\in\N$, $\lambda^{k_n}_i\neq 0$ for some
$i\in\{1,\ldots,k_n\}$.

\no For all $n\in\N$, choose $i_n\in\{1,\ldots,k_n\}$ such that
\beq
|\l^{k_n}_{i_n}|=\max\{|\l^{k_n}_i|\ :\ i=1,\ldots,k_n\}.
\eeq
Thus, we have $\l^{k_n}_{i_n}\neq 0$ $\forall n\in\N$, so
we can choose a sequence $(\hat \e_n)_n$ of positive numbers such
that
\beq
\label{haten}
\lim_{n\to\infty} {\hat \e_n \, k_n\over |\l^{k_n}_{i_n}|}=0
\eeq
(notice that $\hat\e_n\to 0$ as $n\to\infty$ because
$\l^{k_n}_{i_n}\to 0$).
Now, set $\hat\l_n={\l^{k_n}_{i_n}\over |\l^{k_n}_{i_n} |}$ and
  consider the points $\hat y^{k_n}_1,\ldots,\hat y^{k_n}_{k_n}$ in 
$\R^N$ such that 
\beq
\hat y^{k_n}_{i_n}=|x^{k_n}_{i_n}|\, {x^{k_n}_{i_n}+\hat \e_n\hat\l_n\over
    |x^{k_n}_{i_n}+\hat \e_n\hat\l_n |}\ \mbox{ and }\hat
  y^{k_n}_i=x^{k_n}_i\ \mbox{ for }i\neq i_n, \ i\in\{1,\ldots,k_n\}.
\eeq

\no The maximality properties of $r_{k_n}$ and
$(x^{k_n}_1,\ldots,x^{k_n}_{k_n})$ described in Proposition \ref{P2.2}
imply that 
\beq
g^{k_n,\bar\sigma}\left(r_{k_n},{\hat y^{k_n}_1\over |\hat y^{k_n}_1 |},\ldots,
{\hat y^{k_n}_{k_n}\over |\hat y^{k_n}_{k_n} |}\right)\le
g^{k_n,\bar\sigma}\left(r_{k_n},{ x^{k_n}_1\over |x^{k_n}_1 |},\ldots,
{ x^{k_n}_{k_n}\over |x^{k_n}_{k_n} |}\right).
\eeq
Let us choose $(y^{k_n}_1,\ldots,y^{k_n}_{k_n})$ in $D_{k_n}$ such
that
\beq
\left[{1\over {k_n}}\sum_{i=1}^{k_n}|
  y^{k_n}_i|^2\right]^{1/2}=r_{k_n},
\eeq
\beq
(1+2\bar\sigma)^{-1}r_{k_n}\le |
y^{k_n}_i|\le(1+2\bar\sigma)r_{k_n},\quad
{ y^{k_n}_i\over | y^{k_n}_i |}=
{\hat y^{k_n}_i\over |\hat y^{k_n}_i |}\quad\mbox{ for
}i=1,\ldots,k_n
\eeq
and 
\beq
f_{k_n}( y^{k_n}_1,\ldots, y^{k_n}_{k_n} )=
g^{k_n,\bar\sigma}\left(r_{k_n},{\hat y^{k_n}_1\over |\hat y^{k_n}_1|},
\ldots,
{\hat y^{k_n}_{k_n}\over |\hat y^{k_n}_{k_n}|}\right).
\eeq
Then, we have
\beq
f_{k_n}(y^{k_n}_1,\ldots, y^{k_n}_{k_n})\le
f_{k_n}(x^{k_n}_1,\ldots,x^{k_n}_{k_n}),
\eeq
which implies
$E(v_n)\le E(u_{k_n})$ for every $v_n\in
S_{y^{k_n}_1,\ldots,y^{k_n}_{k_n}}$ such that $E(v_n)=f_{k_n}(
y^{k_n}_1,\ldots,$ $y^{k_n}_{k_n})$.
Thus, arguing as in  \cite{CMP,CPS1,CPS2} (in particular as in the
proof of Theorem 2.4 in \cite{CMP}), taking into account the 
assumption (\ref{haten}), Lemma \ref{Lw} and (\ref{form}), we have
\beq
\label{liminf}
\liminf_{n\to \infty}{E(v_n)-E(u_{k_n})\over \hat \e_n\,
  |\l^{k_n}_{i_n}|}\ge
{1\over 2}\,
\int_{\R^N}[w^\d(x)]^2dx\cdot\liminf_{n\to\infty}{1\over
  \hat\e_n} (y_{i_n}^{k_n}-x_{i_n}^{k_n})\cdot\hat\lambda_n.
\eeq
Moreover, since 
\beq
y^{k_n}_{i_n}=|y^{k_n}_{i_n}|\cdot {x^{k_n}_{i_n}+\hat\e_n\hat\l_n\over
  |x^{k_n}_{i_n}+\hat\e_n\hat\l_n |},\qquad 
|y^{k_n}_{i_n}|\ge (1+2\bar\sigma)^{-1}\cdot r_{k_n}\quad
\forall n\in\N, 
\eeq
taking into account that $\lim_{n\to\infty}
{|x^{k_n}_{i_n}|\over r_{k_n}}=1$ because of Lemma
\ref{sigma} and that $x^{k_n}_{i_n}\cdot\hat\lambda_n=0$ $\forall
n\in\N$ because of 
Proposition \ref{P4.2}, we obtain by direct computation
\beq
\liminf_{n\to\infty}{1\over\hat\e_n}(y^{k_n}_{i_n}-x^{k_n}_{i_n})\cdot\hat\lambda_n
=\liminf_{n\to\infty} {|y^{k_n}_{i_n}|\over|x^{k_n}_{i_n}+\hat
  \e\,\hat\lambda_n|}
=\liminf_{n\to\infty} {|y^{k_n}_{i_n}|\over r_{k_n}}
\ge(1+2\bar\sigma)^{-1}
\eeq
which, combined with (\ref{liminf}), implies
\beq
\label{1153}
\liminf_{n\to \infty}{E(v_n)-E(u_{k_n})\over \hat \e_n\,
  |\l^{k_n}_{i_n}|}>0.
\eeq
\no It is clear that (\ref{1153}) is a contradiction because $E(v_n)\le
E(u_{k_n})$ for $n$ large enough, so the proof is complete.

\qed

\no {\underline {\sf Proof of Proposition \ref{P1.2}} \hspace{2mm}}
Proposition \ref{P4.3} guarantees the existence of $\bar k\in \N$ and
of a solution $u_k$, $\forall k\ge\bar k$, having all the properties
reported in  Proposition \ref{P1.2}. 
In fact, (\ref{3}), (\ref{5}) and (\ref{6}) follow from Proposition
\ref{ktoinfty} and (\ref{4}) follows easily from Lemma \ref{sigma}.
Moreover, from (\ref{3}), (\ref{5}), (\ref{6}) and (\ref{B}) we infer
that (\ref{7}) holds and $u_k\to 0$ as $k\to\infty$, uniformly on the
compact subsets of $\R^N$ (because $\lim_{|x|\to\infty}w(x)=0$).

\no Moreover (\ref{w}) implies $\lim_{k\to\infty}\|u_k\|_{H^1(\R^N)}=\infty$.

\no Finally, notice that (\ref{3}), (\ref{5}) and (\ref{6}) imply that
$\liminf_{k\to \infty}E(u_k)\ge k'E_\infty(w)$ $\forall k'\in\N$ so, as
$k'\to\infty$, we obtain $\lim_{k\to\infty}E(u_k)=\infty$.

\qed

\no Theorem \ref{T} is a direct consequence of Proposition \ref{P1.2}.


\begin{rem}{\em
\label{R1}
As a consequence of Lemma \ref{sigma}, the bumps of the
solutions of problem (\ref{.}) are asymptotically distributed near
spheres of dimension $N-1$ when the number of the bumps is large
enough.
If we assume that the potential $a(x)$ satisfies in addition suitable
symmetry assumptions, then we can easily construct infinitely many positive
solutions with bumps distributed near a $d$-dimensional sphere for
every integer $d$ such that $1\le d\le N-1$.
For example, if the potential $a(x)$ satisfies also the symmetry
condition
\beq
\label{5.1}
a(x_1,\ldots,x_{d+1},x_{d+2},\ldots,x_N)=
a(x_1,\ldots,x_{d+1},|x_{d+2}|,\ldots,|x_N|)\quad 
\forall (x_1,\ldots,x_N)\in\R^N,
\eeq
then there exist also solutions with the bumps asymptotically distributed
near spheres of the subspace $\{(x_1,\ldots,x_N)\in\R^N$ : $x_i=0$ for
$i\ge d+2\}$.
In fact, for all $k\in\N$, there exist $\theta_1^k,\ldots,\theta_k^k$
in $S^d=\{x=(x_1,\ldots,x_N)\in\R^N$ : $|x|=1$, $x_i=0$ $\forall i\ge
d+2\}$ and $\rho_1^k,\ldots,\rho_k^k$ in $\R^+$ such that, for $k$
large enough, every minimizing function for $E$ in
$S_{\rho_1^k\theta_1^k,\ldots, \rho_k^k\theta_k^k}$ is a solution
of problem (\ref{.}) and satisfies similar properties as in Proposition
\ref{P1.2} and Lemma \ref{sigma}, in particular the properties
\beq
\label{5.2}
\begin{array}{l}
\vspace{2mm}
\lim\limits_{k\to\infty}\min\{\rho_i^k\ :\ i=1,\ldots,k\}=\infty,
\\ \vspace{2mm}
\lim\limits_{k\to\infty} \min\{|\rho_i^k\theta^k_i-\rho_j^k\theta^k_j|\ :\
i,j=1,\ldots,k,\ i\neq j\}=\infty,
\\
\lim\limits_{k\to\infty}
\frac
{\max\{\rho_i^k\ :\ i=1,\ldots,k\}}
{\min\{\rho_i^k\ :\ i=1,\ldots,k\}}
=1,
\qquad\qquad
\lim\limits_{\e\to 0}\lim\limits_{k\to\infty}{N_k'(\theta,\e)\over
  k\,\e^{d}}\ge c'\quad\forall
\theta\in S^d
\end{array}
\eeq
where $c'$ is a positive constant independent of $\theta$ and
$N'_k(\theta,\e)$ denotes the number of elements of the set
$\{\theta^k_i\in S^d\ :\ i=1,\ldots,k,\ \theta^k_i\in B(\theta,\e)\}$.
}
\end{rem}

\no When $d=1$ and we assume in addition that $a(x)$ has radial symmetry in
the variables $x_1,x_2$, then there exist also other solutions
(corresponding to higher critical levels of the energy functional
$E$).
In fact, for all $k\in\N$ and $\rho>0$, consider the $k$ points in
$\R^N$
\beq
\label{5.3}
x_i^{k,\rho}=\left(\rho\,\cos{2\pi i\over k},\rho\,\sin{2\pi i\over
      k},0,\ldots,0\right)\quad\mbox{ for }i=1,\ldots,k.
\eeq
If $( x_1^{k,\rho},\ldots, x_k^{k,\rho})\in D_k$, set
$\vi_k(\rho)=f_k(x_1^{k,\rho},\ldots, x_k^{k,\rho})$.
Then, as in the proof of Theorem \ref{T}, one can show that there
exists $\tilde k\in\N$ such that, for all $k\ge\tilde k$,
$(x_1^{k,\rho},\ldots, x_k^{k,\rho})$ is in the interior of $D_k$ and
there exists $\tilde r_k>0$ such that 
\beq
\vi_k(\tilde r_k)\ge \vi_k(\rho)\qquad\forall\rho>0\ \mbox{ such
  that }( x_1^{k,\rho},\ldots, x_k^{k,\rho})\in D_k.
\eeq
Let $\tilde u^{k,\tilde r_k}$ be any function in $S_{x_1^{k,\tilde
    r_k},\ldots, x_k^{k,\tilde r_k} }$ such that $E(\tilde u^{k,\tilde
  r_k})=f^k( x_1^{k,\tilde r_k},\ldots, x_k^{k,\tilde r_k} )$.
Because of the radial symmetry, it follows that there exists a Lagrange
multiplier $\tilde \mu_k$ such that $\lambda_i^k=\tilde\mu_k\cdot
x_i^{k,\tilde r_k}$ for $i=1,\ldots,k$. 
Thus, arguing by contradiction as in the proof of Proposition
\ref{P4.3}, one can prove that $\tilde \mu_k=0$ for $k$ large enough,
namely $\tilde u^{k,\tilde r_k} $ is a solution of problem (\ref{.}).

\begin{rem}{\em
\label{R2}
Unlike the results proved in \cite{CPS1,CPS2}, Theorem \ref{T} does
not require $\sup_{x\in\R^N}$ $|a(x)-a_\infty|_{L^{N/2}(B(x,1))}$ to be
small and, indeed, it may be arbitrarily large. 
For example, if $\Omega$ is a bounded domain of $\R^N$ and
$a_s(x)=s\bar a(x)+a(x)$ $\forall x\in\R^N$, where $a(x)$ is as in
Theorem \ref{T} and $\bar a(x)$ is a nonnegative function which is
positive only in $\Omega$, then for $k$ large enough and all $s\ge 0$ there
exists a $k$-bump solution $u_{k,s}$; moreover, as $s\to\infty$,
$u_{k,s}$ converges to a $k$-bump solution $\tilde u_k$ in the
exterior domain $\widetilde\Omega=\R^N\setminus\overline\Omega$, with zero
Dirichlet boundary condition (on the other hand, the solution $\tilde
u_k$ may be also obtained directly since our method can be
adapted to deal with Dirichlet problems in exterior domains).
}\end{rem}

\begin{rem}{\em
The method developed in this paper may be also used to construct
sequences $(\hat u_n)_n$ of positive solutions of problem (\ref{.})
which converge in $H^1_{\loc} (\R^N)$ to a positive solution $\hat u$
having infinitely many bumps (while the sequence $(u_k)_{k\geq\bar k}$ given by
Theorem \ref{T} converges to the trivial solution $u\equiv 0$).
The bumps are distributed near infinitely many spheres with center in
the origin.
Since the radius of these spheres may be chosen in infinitely many ways,
we obtain infinitely many positive solutions having
infinitely many positive bumps (while the result presented in
\cite{CPS2} guarantees only the existence of one solution having this
property, under the additional assumption that
$\sup_{x\in\R^N}|a(x)-a_\infty|_{L^{N/2}(B(x,1))}$ is small enough).
}\end{rem}





\vspace{4mm}

{\Large{\bf{Appendix} }}

\vspace{4mm}

\no Here we report a table of the main notations used in the paper and we
indicate the formula where every notation is used the first time.

\begin{longtable}{c|c||c||c|c||c||c|c}
\endfirsthead
\multicolumn{8}{r}{\textit{(Continue to next page)}} 
\endfoot
\multicolumn{8}{l}{\textit{(Continue from previous page)}} 
\endhead
\multicolumn{8}{c}{}\\
\endlastfoot
$a(x)$             &             (\ref{.})           & &
$a_\infty$          &            (\ref{1})           & &
$E$                &            (\ref{E})          \\
\hline
$Z$             &            (\ref{2.2})           & &
$R_\d$          &            (\ref{2.3})          & &
$D_k$                &            (\ref{2.3})          \\
\hline
$\beta_i$                &            (\ref{2.4})          & & 
{\small{$S_{x_1,\ldots,x_k}$ }}          &            (\ref{2.6})     & & 
$f_k$          &            (\ref{df})            \\   
\hline
$\l_i$                &            (\ref{lambda})          & &
{\small{$D^{k,\sigma}(\rho,\theta_1,\ldots,\theta_k)$}} &       (\ref{2.17})   & &        
$D^{k,\sigma}$             &            (\ref{2.18})          \\
\hline
$g^{k,\sigma}$          &            (\ref{g})           & &
$r_k$           &            (\ref{2.20})           & &
$x^k_i$ &       (\ref{2.20})         \\
\hline
$E_\infty$             &            (\ref{kEinfty})           & &
$\bar\sigma_n$          &            (\ref{3.6})           & & 
$r_{k,\bar\sigma_n}$   &            (\ref{3.6})          \\
\hline
$x^{k,\bar\sigma_n}_i$ &       (\ref{3.6})          & &
$\sigma_n$             &            (\ref{3.10})       & &     
$V_n$          &            (\ref{Vn})          \\
\hline
$S(x,r)$           &            (\ref{Vn})          & &
$\nu_n$ &       (\ref{gamma})             & &
$\gamma_n$             &            (\ref{gamma})        \\
\hline
$N_{k_n}(x,\e)$          &            (\ref{e3.15})           & &
$\mu_k$           &            (\ref{lambdamu})        & &   
$\bar\sigma$ &       (\ref{<})         \\
\hline
$w_{k,i} $    &       (\ref{muk})          & &
$C_{n,i} $    &       (\ref{4.15})          & & 
$x^{k,\rho}_i$ &        (\ref{5.3})         \\

\end{longtable}

\vspace{4mm}

{\bf Acknowledgements}.
This work is supported by  ``Gruppo
Nazionale per l'Analisi Matematica, la Probabilit\`a e le loro
Applicazioni (GNAMPA)'' of the {\em Istituto Nazionale di Alta Matematica
(INdAM)}.


\section*{Compliance with Ethical Standards}

We wish to confirm that there are no known conflicts of interest
associated with this publication and there has been no significant
financial support for this work that could have influenced its
outcome.



\end{document}